\documentclass[10pt]{amsart}

\newtheorem{cotmb}{Corollary}
\newtheorem{prop}{Proposition}
\newcommand{\cal}{\mathcal}

\newtheorem{tw}{Theorem}
\newtheorem{lemma}{Lemma}

\usepackage[mathscr]{euscript}
\usepackage[normalem]{ulem}

\vspace{0.5cm}

\title{Harmonic almost Hermitian structures}

\author{Johann Davidov}
\thanks{The author is partially supported by  the National Science
Fund, Ministry of Education and Science of Bulgaria under contract
DFNI-I 02/14. }

\address{Johann Davidov\\Institute of Mathematics and Informatics \\
Bulgarian Academy of Sciences\\ Acad. G.Bonchev Str. Bl.8\\ 1113
Sofia\\ Bulgaria\\  \newline \centerline{and} \newline
University of Structural Engineering and Architecture "L.Karavelov",
175 Suhodolska St., 1373 Sofia, Bulgaria }

\begin{document}

\begin{abstract}
This is a survey of old and new results on the problem when a
compatible almost complex structure on a Riemannian manifold is a
harmonic section or a harmonic map from the manifold into its
twistor space. In this context, special attention is paid to the
Atiyah-Hitchin-Singer and Eells-Salamon almost complex structures on
the twistor space of an oriented Riemannian four-manifold.

\vspace{0,1cm} \noindent 2010 {\it Mathematics Subject
Classification}. Primary 53C43, Secondary 58E20, 53C28

\vspace{0,1cm} \noindent {\it Key words: almost complex structures,
twistor spaces, harmonic maps}
\end{abstract}

\thispagestyle{empty}

\maketitle
\vspace{0.5cm}

\section{Introduction}

Recall that an almost complex structure on a Riemannian manifold
$(N,h)$ is called almost Hermitian (or, compatible) if it is
$h$-orthogonal. If a Riemannian manifold admits an almost Hermitian
structure, it possesses many such structures (cf. Section~\ref{SV}).
Thus,  it is natural to look for "reasonable" criteria that
distinguish some of these structures. A natural way to obtain such
criteria is to consider the almost Hermitian structures on $(N,h)$
as sections of its twistor bundle $\pi:{\mathscr T}\to N$ whose
fibre at a point $p\in N$ consists of all $h$-orthogonal complex
structures $J_p:T_pN\to T_pN$  on the tangent space of $N$ at $p$.
The fibre of the bundle ${\mathscr T}$ is the compact Hermitian
symmetric space $O(2m)/U(m)$, $2m=dim\,N$, and its standard metric
$-\frac{1}{2}Trace\,J_1\circ J_2$ is K\" ahler-Einstein. The twistor
space ${\mathscr T}$ admits a natural Riemannian metric $h_1$ such
that the projection map $\pi:({\mathscr T},h_1)\to (N,h)$ is a
Riemannian submersion with totally geodesic fibres. This metric is
compatible with the natural almost complex structures on ${\mathscr
T}$, which  have been introduced  by Atiyah-Hitchin-Singer
\cite{AHS} and Eells-Salamon \cite{ES} in the case $dim\,N=4$.

If $N$ is oriented, the twistor space ${\mathscr T}$ has two
connected components often called positive and negative twistor
spaces of $(N,h)$; their sections are almost Hermitian structures
yielding the orientation and, respectively, the opposite orientation
of $N$.

    E. Calabi and H. Gluck \cite{CG} have proposed to single out those almost Hermitian structures
$J$ on $(N,h)$, whose image $J(N)$ in ${\mathscr T}$ is of minimal volume with respect to the
 metric $h_1$. They have proved that the standard
almost Hermitian structure on the $6$-sphere $S^6$, defined by means
of the Cayley numbers, can be characterized by that property.

    Motivated by harmonic map theory, C. M. Wood \cite{W1,W2} has suggested to consider as "optimal"
those almost-Hermitian structures $J:(N,h)\to ({\mathscr T},h_1)$,
which are critical points of the energy functional under variations
through sections of ${\mathscr T}$. In general, these critical
points are not harmonic maps,  but, by analogy, in \cite{W2} they
are referred to as "harmonic almost complex structures". They are
also called "harmonic sections" \cite{W1}, a term, which is more
appropriate in the context of this article.

Forgetting the bundle structure of ${\mathscr T}$, we can  consider
almost Hermitian structures that are critical points of the energy
functional under variations through all maps $N\to{\mathscr T}$.
These structures are genuine harmonic maps from $(N,h)$ into
$({\mathscr T},h_1)$ and we refer to \cite{EL} for basic fact about
such maps.

The main goal of this paper is to survey results about the
harmonicity (in both senses) of the Atiyah-Hitchin-Singer and
Eells-Salamon almost Hermitian structures on the twistor space of an
oriented four-dimensional Riemannian manifold, as well as almost
Hermitian structures on such a manifold.

In Section 2 we recall some basic facts about the twistor spaces of
even-dimensional Riemannian manifolds. Special attention is paid to
the twistor spaces of oriented four-dimensional manifolds. In
Sections 3 and 4 we discuss the energy functional on sections of a
twistor space, i.e. almost Hermitian structures on the base
Riemannian manifold. We state the Euler-Lagrange equation for such a
structure to be a critical point of the energy functional (a
harmonic section) obtained by Wood \cite{W1,W2}. Several examples of
non-K\"ahler almost Hermitian structures, which are harmonic
sections are given. K\"ahler structures are absolute minima of the
energy functional. G. Bor, L. Hern\'andez-Lamoneda and M. Salvai
\cite{BLS} have given sufficient conditions for an almost Hermitian
structure to be a minimizer of the energy functional. Their result
(in fact, part of it) is presented in Section 4 and is used to
supply examples of non-K\"ahler minimizers based on works by C.
LeBrun\cite{LeBr} and I. Kim \cite{Kim}. Section 4 ends with a lemma
from \cite{DM02}, which rephrases the Euler-Lagrange equation for an
almost Hermitian structure $(h,J)$ on a manifold $N$ in terms of the
fundamental $2$-form of $(h,J)$ and the curvature of $(N,h)$. This
lemma has been used in \cite{DM02} to show that the
Atiyah-Hitchin-Singer almost Hermitian structure ${\mathcal J}_1$ on
the negative twistor space ${\mathcal Z}$ of an oriented Riemannian
$4$-manifold $(M,g)$ is a harmonic section if and only if the base
manifold $(M,g)$ is self-dual, while the Eells-Salamon structure
${\mathcal J}_2$ is  a harmonic section if and only $(M,g)$ is
self-dual and of constant scalar curvature. The main part of the
proof of this result (slightly different from the proof in
\cite{DM02}) is presented in Section 5. In this context, it is
natural to ask when ${\mathcal J}_1$ and ${\mathcal J}_2$ are
harmonic maps into the twistor space of ${\cal Z}$. Recall that a
map between Riemannian manifolds is harmonic if and only if the
trace of its second fundamental form vanishes. Section 6 contains a
computation of the second fundamental form of a map from a
Riemannian manifold $(N,h)$ into its twistor space $({\mathscr
T},h_1)$. The corresponding formula obtained in this section is used
in Section 7 to give an answer to the question above: ${\mathcal
J}_1$ or ${\mathcal J}_2$ is a harmonic map if and only if  $(M,g)$
is either self-dual and  Einstein, or locally is the product of an
open interval on ${\mathbb R}$ and a $3$-dimensional Riemannian
manifold of constant curvature. A sketch of  the proof, involving
the main theorem of Section 5 and several technical lemmas, is given
in Section 7 following \cite{DM16}. In Section 8 we give  geometric
conditions on a four-dimensional almost Hermitian manifold $(M,g,J)$
under which the almost complex structure $J$ is a harmonic map of
$(M,g)$ into the positive twistor space $({\mathcal Z}_{+},h_1)$,
$M$ being considered with the orientation induced by $J$. We also
find conditions for minimality of the submanifold $J(M)$ of the
twistor space ${\mathcal Z}_{+}$. As is well-known, in dimension
four, there are three basic classes in the Gray-Hervella
classification \cite{GH} of almost Hermitian structures:  Hermitian,
almost K\"ahler (symplectic) and K\"ahler structures. As for a
manifold of an arbitrary dimension, if $(g,J)$ is K\"ahler, the map
$J:(M,g)\to ({\mathcal Z}_{+},h_1)$ is a totally geodesic isometric
imbedding.  In the case of a Hermitian structure, we express the
conditions for harmonicity or minimality of $J$ in terms of the Lee
form, the Ricci and star-Ricci tensors of $(M,g,J)$, while for an
almost K\"ahler structure the conditions are in terms of the Ricci,
star-Ricci and Nijenhuis tensors. Several examples illustrating
these results are discussed at the end of Section 8, among them a
Hermitian structure that is a harmonic section of the twistor bundle
${\cal Z}_{+}$ and a minimal isometric imbedding in it, but not a
harmonic map.

\smallskip

\noindent {\bf Acknowledgment}. This paper is an expanded version on
the author's talk at the INdAM workshop in Rome, November, 16-20,
2015, on the occasion of the sixtieth birthday of Simon Salamon. The
author would like to express his gratitude to the organizers and
Simon for the invitation to take part in the workshop and for the
wonderful and stimulating environment surrounding it. Special thanks are also due to
the referee whose remarks helped to improve the final version of the article.

\smallskip

\section{The twistor space of an even-dimensional Riemannian manifold}
Denote by $F({\mathbb R^{2m}})$ the set of complex structures on
${\mathbb R^{2m}}$ compatible with its standard metric $g$. This set
has the structure of an imbedded  submanifold of the vector space
$so(2m)$ of skew-symmetric endomorphisms of ${\mathbb R^{2m}}$. The
tangent space  of $F({\mathbb R^{2m}})$ at a point $J$ consists of
all skew-symmetric endomorphisms of ${\mathbb R^{2m}}$
anti-commuting with $J$. Then we can define an almost complex
structure on the manifold $F({\mathbb R^{2m}})$ setting
$$
{\mathcal J}Q=JQ \; ~~\mbox {for} ~~\; Q\in T_JF({\mathbb R^{2m}}).
$$
This almost complex structure is compatible with the standard metric
$$G(A,B)=-\frac{1}{2m}Trace\,AB$$ of $so(2m)$, where the factor $1/2m$
is chosen so that every $J\in F({\mathbb R^{2m}})$ should have unit
norm. In fact, the almost Hermitian structure $(G,{\cal J})$ is
K\"ahler-Einstein.

 The group $O(2m)$ acts on $F({\mathbb R^{2m}})$ by
conjugation and the isotropy subgroup at the standard complex
structure $J_0$ of ${\mathbb R^{2m}}\cong {\mathbb C^{m}}$ is
$U(m)$. Therefore $F({\mathbb R^{2m}})$ can be identified with the
homogeneous space $O(2m)/U(m)$.

Note also that the manifold $F({\mathbb R^{2m}})$ has two connected
components: if we fix an orientation on ${\mathbb R}^ {2m}$, these
components consists of all complex structures on ${\mathbb R}^{2m}$
compatible with the metric $g$ and inducing $\pm$ the orientation of
${\mathbb R}^{2m}$; each of them has the homogeneous representation
$SO(2m)/U(m)$.

\smallskip

The twistor space of an even-dimensional Riemannian manifold
$(N,h)$, $dim\,N=2m$, is the bundle $\pi:{\mathscr T}\to N$,  whose
fibre at every point $p\in N$ is the space of compatible complex
structures on the Euclidean vector space $(T_pN,h_p)$. This is the
associated bundle
$$
{\mathscr T}=O(N)\times_{O(2m)}F({\mathbb R^{2m}})
$$
where $O(N)$ is the principal bundle of orthonormal frames on $N$.
Since the bundle $\pi:{\mathscr T}\to N$ is associated to $O(N)$,
the Levi-Civita connection on $(N,h)$ gives rise to a splitting
${\mathcal V}\oplus {\mathcal H}$ of the tangent bundle of the
manifold ${\mathscr T}$. Using this splitting, we can define a
natural $1$-parameter family of Riemannian metrics $h_t$, $t>0$, as
follows. For every $J\in {\mathscr T}$, the horizontal subspace
${\mathcal H}_J$ of $T_J{\mathscr T}$ is isomorphic  to the tangent
space $T_{\pi(J)}N$ via the differential $\pi_{\ast J}$ and the
metric $h_t$ on ${\cal H}_J$ is the lift of the metric $h$ on
$T_{\pi(J)}N$, $h_t|{\mathcal H}_J=\pi^{\ast} h$. The vertical
subspace ${\cal V}_J$ of $T_J{\cal T}$ is the tangent space at $J$
to the fibre through $J$ of the bundle ${\cal T}$  and $h_t|{\cal
V}_J$ is defined as $t$ times the metric $G$ of this fibre. Finally,
the horizontal space ${\cal H}_J$ and the vertical space ${\cal
V}_J$ are declared to be orthogonal. Then, by  the Vilms theorem
\cite{V}, the  projection $\pi:({\mathscr T},h_t)\to (N,h)$ is a
Riemannian submersion with totally geodesic fibres (this can also be
proved directly).

It is often convenient to consider  ${\mathscr T}$ as a submanifold
of the bundle
$$
\pi:A(TN)=O(N)\times_{O(2m)}so(2m)\to N
$$
of skew-symmetric endomorphisms of $TN$. The inclusion of ${\mathscr
T}$ into $A(TN)$ is fibre-preserving and for every $J\in{\mathscr
T}$ the horizontal subspace ${\cal H}_J$ of $T_J{\mathscr T}$
coincides with the horizontal subspace of $T_JA(TN)$ with respect to
the connection induced by the Levi-Civita connection of $(N,h)$
since the inclusion of $F({\mathbb R^{2m}})$ into $so(2m)$ is
$O(2m)$-equivariant; the vertical subspace ${\cal V}_J$ of
$T_J{\mathscr T}$ is the subspace of the fibre $A(T_{\pi(J)}N)$ of
$A(TN)$ through $J$ consisting of the skew-symmetric endomorphisms
of $T_{\pi(J)}N$ anti-commuting with $J$.

If the manifold $N$ is oriented, its twistor space has two connected
components,  the spaces of compatible complex structures on tangent
spaces of $N$ yielding the given, or the opposite orientation of
$N$. These are often called the positive, respectively, the negative
twistor space.

\subsection {The twistor space of an oriented four-dimensional
Riemannian manifold} In dimension four, each of the two connected
components $F({\mathbb R}^4)$ can be identified with the unit sphere
$S^2$. It is often convenient to describe this identification in
terms of the space $\Lambda^2{\mathbb R}^4$. The metric $g$ of
${\mathbb R}^4$ induces a metric on $\Lambda^2{\mathbb R}^4$ given
by
\begin{equation}\label{g^2}
g(x_1\wedge x_2,x_3\wedge
x_4)=\frac{1}{2}[g(x_1,x_3)g(x_2,x_4)-g(x_1,x_4)g(x_2,x_3)],
\end{equation}
the  factor $1/2$ being chosen in consistency with \cite{DM91,DM02}.
Consider the isomorphisms $so(4)\cong \Lambda^2{\mathbb R}^4$
sending $\varphi\in so(4)$ to the $2$-vector $\varphi^{\wedge}$ for
which
$$
2g(\varphi^{\wedge},x\wedge y)=g(\varphi x,y),\quad x,y\in {\mathbb
R}^4.
$$
This isomorphism is an isometry with respect to the metric $G$ on
$so(4)$ and the metric $g$ on $\Lambda^2{\mathbb R}^4$. Given $a\in
\Lambda^2{\mathbb R}^4$, the skew-symmetric endomorphism of
${\mathbb R}^4$ corresponding to $a$ under the inverse isomorphism
will be denoted by $K_a$.

Fix an orientation on ${\mathbb R}^4$ and denote by
$F_{\pm}({\mathbb R}^4)$ the set of complex structures on ${\mathbb
R}^4$ compatible with the metric $g$ and inducing $\pm$ the
orientation of ${\mathbb R}^4$.

The Hodge star operator defines an endomorphism $\ast$ of
$\Lambda^2{\mathbb R}^4$ with $\ast^2=Id$. Hence we have the
orthogonal decomposition
$$
\Lambda^2{\mathbb R}^4=\Lambda^2_{-}{\mathbb
R}^4\oplus\Lambda^2_{+}{\mathbb R}^4,
$$
where $\Lambda^2_{\pm}{\mathbb R}^4$ are the subspaces of
$\Lambda^2{\mathbb R}^4$ corresponding to the $(\pm 1)$-eigenvalues
of the operator $\ast$. Let $(e_1,e_2,e_3,e_4)$ be an oriented
orthonormal basis of ${\mathbb R}^4$. Set
\begin{equation}\label{s-basis}
s_1^{\pm}=e_1\wedge e_2 \pm e_3\wedge e_4, \quad s_2^{\pm}=e_1\wedge
e_3\pm e_4\wedge e_2, \quad s_3^{\pm}=e_1\wedge e_4\pm e_2\wedge
e_3.
\end{equation}
Then $(s_1^{\pm},s_2^{\pm},s_3^{\pm})$ is an orthonormal basis of
$\Lambda^2_{\pm}{\mathbb R}^4$.

It is easy to see that the isomorphism $\varphi\to \varphi^{\wedge}$
identifies $F_{\pm}({\mathbb R}^4)$ with the unit sphere
$S(\Lambda^2_{\pm}{\mathbb R}^4)$ of the Euclidean vector space
$(\Lambda^2_{\pm}{\mathbb R}^n,g)$. Under this isomorphism, if $J\in
F_{\pm}({\mathbb R}^4)$, the tangent space $T_JF({\mathbb
R}^4)=T_JF_{\pm}({\mathbb R}^4)$ is identified with the orthogonal
complement $({\mathbb R} J)^{\perp}$ of the space ${\mathbb R}J$ in
$\Lambda^2_{\pm}{\mathbb R}^4$.

\begin{lemma}
The orientation on $\Lambda^2_{\pm}{\mathbb R}^4 $ determined by the
basis $s_1^{\pm},s_2^{\pm},s_3^{\pm}$ defined by means of an
oriented orthonormal basis $\{e_1,...,e_4\}$ of ${\mathbb R}^4$ does not
depend on the choice of $\{e_1,...,e_4\}$.
\end{lemma}

\noindent {\bf Proof}. Let $\{s_i'=s_i'^{+}\}$ and $\{s_i=s_i^{+}\}$
be the bases of $\Lambda^2_{+}{\mathbb R}^4$ defined by means of two
oriented orthonormal bases $\{e_1',...,e_4'\}$ and $\{e_1,...,e_4\}$ of
${\mathbb R}^4$. Denote by $A\in SO(4)$ the transition matrix
between these bases. Thanks to L. van Elfrikhof (1897), it is well-known
that every matrix $A$ in $SO(4)$ can be represented as the product
$A=A_1A_2$ of two $SO(4)$-matrices of the following types
\begin{equation}\label{isoclinic}
A_1=\left(\begin{array}{cccc}
                   a & -b & -c & -d \\
                   b & a & -d & c \\
                   c & d & a & -b \\
                   d & -c & b & a \\
                  \end{array} \right),\quad
A_2=\left(\begin{array}{cccc}
                   p & -q & -r & -s \\
                   q & p & s & -r \\
                   r & -s & p & q \\
                   s & r & -q & p \\
                  \end{array} \right),
\end{equation}
where $a,...,d,p,...,s$ are real numbers with $a^2+b^2+c^2+d^2=1$,
$p^2+q^2+r^2+s^2=1$ (isoclinic representation). For an endomorphism
$L$ of ${\mathbb R}^4$, denote by $\Lambda_L$ the induced
endomorphism on $\Lambda^2{\mathbb R}^4$ defined by
$\Lambda_L(X\wedge Y)=L(X)\wedge L(Y)$. Denote again by $A$ the
isomorphism of ${\mathbb R}^4$ with matrix $A$ with respect to the
basis $e_1,...,e_4$. Then $s_i'=\Lambda_A(s_i)=\Lambda_{A_1}\circ
\Lambda_{A_2}(s_i)$, $i=1,2,3$. One easily computes that
$\Lambda_{A_2}(s_i)=s_i$, and that $\Lambda_{A_1}(s_i)$ is a basis
of $\Lambda^2_{+}{\mathbb R}^4$ whose transition matrix is
$$
\left( \begin{array}{cccc}
                                       a^2+b^2-(c^2+d^2) & -2(ad-bc) & 2(ac+bd)  \\
                                       2(ad+bc) & a^2+c^2-(b^2+d^2) & -2(ab-cd) \\
                                       -2(ac-bd) & 2(ab+cd) & a^2+d^2-(b^2+c^2)\\
                                                                            \end{array}
                                   \right).
$$
The determinant of the latter matrix is $(a^2+b^2+c^2+d^2)^3=1$.
This proves the statement for $\Lambda^2_{+}{\mathbb R}^4$. Changing
the orientation of ${\mathbb R}^4$ interchanges the roles of
$\Lambda^2_{+}{\mathbb R}^4$ and $\Lambda^2_{-}{\mathbb R}^4$.
Therefore, the statement holds  for $\Lambda^2_{-}{\mathbb R}^4$ as
well.

The orientation of $\Lambda^2_{\pm}$ induced by a basis
$\{s_1^{\pm},s_2^{\pm},s_3^{\pm}\}$ will be called  "canonical".

\smallskip

\noindent {\bf Remark}. The map  assigning the coset of the matrix  above in $SO(3)/SO(2)=S^2$ to the unit quaternion  $q=a+ib+jc+kd$ is the Hopf map $S^3\to S^2$.

\smallskip

Consider the $3$-dimensional Euclidean space
$(\Lambda^2_{\pm}{\mathbb R}^4,g)$ with its canonical orientation
and denote by $\times$ the usual vector-cross product on it. Then,
if $a,b\in\Lambda^2_{\pm}{\mathbb R}^4$, the isomorphism
$\Lambda^2{\mathbb R}^4\cong so(4)$ sends $a\times b$ to
$\pm\frac{1}{2}[K_a,K_b]$. Thus, if $J\in F_{\pm}({\mathbb R}^4)$
and $Q\in T_JF({\mathbb R}^4)=T_JF_{\pm}({\mathbb R}^4)$, we have
\begin{equation}\label{calJ}
({\mathcal J}Q)^{\wedge}=\pm ( J^{\wedge}\times  Q^{\wedge}).
\end{equation}

\smallskip

Now let $(M,g)$ be an oriented Riemannian manifold of dimension
four.

According to the considerations above, the twistor space of such a
manifold has two connected components, which can be identified with
the unit sphere subbundles ${\cal Z}_{\pm}$ of the bundles
$\Lambda^2_{\pm}TM\to M$,  the eigensubbundles of the bundle
$\pi:\Lambda^2TM\to M$ corresponding to the eigenvalues $\pm 1$ of
the Hodge star operator. These are the positive and the negarive
twistor spaces of $(M,g)$. If $\sigma\in{\mathcal Z}_{\pm}$, then
$K_{\sigma}$ is a complex structure on the vector space
$T_{\pi(\sigma)}M$ compatible with the metric $g$ and $\pm$ the
orientation of $M$. Note that, since
$g(K_{\sigma}X,Y)=2g(\sigma,X\land Y)$, the $2$-vector $2\sigma$ is
dual  to  the fundamental $2$-form of $(g,K_{\sigma })$.

\smallskip

The Levi-Civita connection $\nabla$ of $M$ preserves the bundles
$\Lambda^2_{\pm}TM$, so it induces a metric connection on each of
them denoted again by $\nabla$. The  horizontal distribution of
$\Lambda^2_{\pm}TM$ with respect to $\nabla$ is tangent to the
twistor space ${\cal Z}_{\pm}$.

The manifold ${\cal Z}_{\pm}$ admits two almost complex structures
${\cal J}_1$ and ${\cal J}_2$ compatible with each metric $h_t$
defined, respectively, by Atiyah-Hitchin-Singer \cite{AHS}, and
Eells-Salamon \cite{ES}.  On a vertical space ${\cal V}_J$, ${\cal
J}_1$ is defined to be the complex structure ${\cal J}_J$ of the
fibre through $J$, while ${\cal J}_2$ is defined as the conjugate
complex structure, i.e. ${\cal J}_{2}|{\cal V}_J=-{\cal J}_J$. On a
horizontal space ${\cal H}_J$, ${\cal J}_1$ and ${\cal J}_2$ are
both defined to be the lift to ${\cal H}_J$ of the endomorphism $J$
of $T_{\pi(J)}M$. Thus, if $\sigma\in{\cal Z}_{\pm}$
$${\cal J}_{k}|{\mathcal H}_{\sigma}=(\pi_{\ast}|{\mathcal H}_{\sigma})^{-1}\circ K_{\sigma}
\circ\pi_{\ast}|{\mathcal H}_{\sigma}.
$$
$$
{\cal J}_{k}V=(-1)^{k}\sigma\times V \;~~\mbox {for} ~~\; V\in
\mathcal{V_{\sigma}},\;\;\;\;k=1,2.
$$

\smallskip

Let $R$ be the curvature tensor of the Levi-Civita connection of
$(M,g)$; we adopt the following definition for the curvature tensor
$R$: $R(X,Y)=\nabla_{[X,Y]}-[\nabla_{X},\nabla_{Y}]$. Then the
curvature operator ${\cal R}$ is the self-adjoint endomorphism of
$\Lambda ^{2}TM$ defined by
$$
 g({\cal R}(X\land Y),Z\land T)=g(R(X,Y)Z,T),\quad X,Y,Z,T\in TM.
$$
The curvature tensor of the connection on the bundle $\Lambda^2TM$ induced by the Levi-Civita connection  on $TM$ will also
be denoted by $R$.

\smallskip

The following easily verified formulas are useful in various computations on ${\mathcal Z}_{\pm}$.
\begin{equation}\label{eq 2.5}
g(R(a)b,c)=\pm g({\cal R}(a),b\times c))
\end{equation}
for $a\in\Lambda^2T_{p}M$, $b,c\in \Lambda ^{2}_{\pm}T_{p}M$,
\begin{equation}\label{KK}
K_{b}\circ K_{c}=-g(b,c)Id\pm K_{b\times c},\quad b,c\in \Lambda ^{2}_{\pm}T_{p}M.
\end{equation}
\begin{equation}\label{eq 2.6}
g(\sigma\times V,X\land K_{\sigma }Y)=%
g(\sigma\times V,K_{\sigma }X\land Y)=\pm g(V,X\land Y)
\end{equation}
for  $\sigma\in{\mathcal Z}_{\pm}$, $V\in {\cal V}_{\sigma}$, $X,Y\in T_{\pi(\sigma)}M$.

\smallskip

Denote by ${\cal B}:\Lambda^2TM\to \Lambda^2TM$ the endomorphism
corresponding to the traceless Ricci tensor. If  $s$ denotes the
scalar curvature of $(M,g)$ and $\rho:TM\to TM$ is the Ricci
operator, $g(\rho(X),Y)=Ricci(X,Y)$, we have
$$
{\mathcal B}(X\wedge Y)=\rho(X)\wedge
Y+X\wedge\rho(Y)-\frac{s}{2}X\wedge Y.
$$
Note that ${\mathcal B}$ sends  $\Lambda^2_{\pm}TM$ into
$\Lambda^2_{\mp}TM$. Let ${\cal W}: \Lambda^2TM\to \Lambda^2TM$ be
the endomorphism corresponding to the Weyl conformal tensor. Denote
the restriction of ${\cal W}$ to $\Lambda^2_{\pm}TM$ by ${\cal
W}_{\pm}$, so ${\cal W}_{\pm}$ sends $\Lambda^2_{\pm}TM$ to
$\Lambda^2_{\pm}TM$ and vanishes on $\Lambda^2_{\mp}TM$.

 It is well known that the curvature operator decomposes as (\cite{ST}, see
e.g. \cite[Chapter 1 H]{Besse})
\begin{equation}\label{dec}
{\cal R}=\frac{s}{6}Id+{\cal B}+{\cal W}_{+}+{\cal W}_{-}
\end{equation}
Note that this differ by  a
factor $1/2$ from \cite{Besse} because of the factor $1/2$ in our
definition of the induced metric on $\Lambda^2TM$.

The Riemannian manifold $(M,g)$ is Einstein exactly when ${\cal
B}=0$. It is called self-dual (anti-self-dual), if ${\cal W}_{-}=0$
(resp. ${\cal W}_{+}=0$).

It is a well-known result by Atiyah-Hitchin-Singer  \cite{AHS} that
the almost complex structure ${\cal J}_1$ on ${\cal Z}_{-}$ (resp.
${\cal Z}_{+}$) is integrable (i.e. comes from a complex structure)
if and only if $(M,g)$ is self-dual (resp. anti-self-dual). On the
other hand the almost complex structure ${\cal J}_2$ is never
integrable by a result of Eells-Salamon \cite{ES}, but nevertheless
it plays a useful role in harmonic map theory.

\section {The standard variation with compact support of an almost
Hermitian structure through sections of the twistor space}\label{SV}

Now suppose that $(N,h)$ is a Riemannian manifold, which admits an
almost Hermitian structure $J$, i.e. a section of the bundle
$\pi:{\mathscr T}\to N$. Take a section $V$ with compact support $K$
of the bundle $J^{\ast}{\cal V}\to N$, the pull-back under $J$ of
the vertical bundle ${\cal V}\to {\mathscr T}$. There exists  an
$\varepsilon>0$ such that,  for every point $I$ of the compact set
$J(K)$,  the exponential map $exp_{I}$ is a diffeomorphism of the
$\varepsilon$-ball in $T_I{\mathscr T}$. The function $||V||_{h_1}$
is bounded on $N$, so there exists a number  $\varepsilon'>0$ such
that $||sV(p)||_{h_1}<\varepsilon$ for every $p\in N$ and
$s\in(-\varepsilon',\varepsilon')$. Set $J_s(p)=exp_{J(p)}[sV(p)]$
for $p\in N$ and $s\in(-\varepsilon',\varepsilon')$. For every fixed
$p\in N$, the curve $s\to exp_{J(p)}[sV(p)]$ is a geodesic with
initial velocity vector $V(p)$ which is tangent to the fibre
${\mathscr T}_p$ of ${\mathscr T}$ through $J(p)$. Since this fibre
is a totally geodesic submanifold, the whole curve lies in it. Hence
$J_s$ is a section of ${\mathscr  T}$, i.e. an almost Hermitian
structure on $(N,h)$, such that $J_s=J$ on $N\setminus K$.

In particular, this shows that if $(N,h)$ admits a compatible almost
complex structure $J$, then it possesses many such structures.

\section{The energy functional on sections of the twistor space}

If $D$ is a relatively compact open subset of a Riemannian manifold
$(N,h)$, the energy functional assigns to every compatible almost
complex structure $J$ on $(N,h)$, considered as a map $J:(N,h)\to
({\mathscr T},h_t)$, the integral
$$
E_{\Omega}(J)=\int_{D}||J_{\ast}||^2_{h,h_t}vol
$$
where the norm is taken with respect to $h$ and $h_t$.

A compatible almost complex structure $J$ is said to be a {\it
harmonic section} ("a harmonic almost complex structure" in the
terminology of \cite{W2}), if for every $D$ it is a critical point
of the energy functional under variations of $J$ through sections of
the twistor space of $(N,h)$.

We have $J_{\ast}X=\nabla_XJ+X^h$ for every $X\in TN$ where $X^h$ is
the horizontal lift of $X$ (and $\nabla_XJ$ is the vertical part of
$J_{\ast}X)$. Therefore $||J_{\ast}||^2=t||\nabla
J||^2_h+(dim\,N)vol(D)$. It follows that the critical points of
the energy functional $E_{\Omega}$ coincide with the critical points
of the vertical energy functional
$$
J\to \int_{D}||\nabla J||^2vol
$$
and do not depend on the particular choice of the parameter $t$.
Another obvious consequence is that the K\"ahler structures provide
the absolute minimum of the energy functional.

The Euler-Lagrange equation for the critical points of the energy
functional under variations through sections of the twistor bundle
has been obtained by C. Wood.

\begin{tw}\rm{(\cite{W1,W2})}
A compatible almost complex structure $J$ is a harmonic section if
and only if
$$
[J,\nabla^{\ast}\,\nabla J]=0,
$$
where $\nabla^{\ast}$ is the formal adjoint operator of $\nabla$.
\end{tw}

\smallskip

\noindent {\bf Remark}. Suppose that $N$ is oriented and $J$ is an
almost Hermitian structure on $(N,h)$ yielding the orientation of
$N$, so it is a section of the positive twistor bundle ${\mathscr
T}_{+}$. Every variation of $J$ with compact support consisting of
sections of the total twistor space ${\mathscr T}$ contains a
subvariation consisting of sections of ${\mathscr T}_{+}$. Thus  $J$
is a critical point of the energy functional under variations
 with compact support through sections of the total twistor space ${\mathscr T}$,
 if and only if it is a critical point under variations through sections of
${\mathscr T}_{+}$.

\smallskip

\noindent {\bf Examples} of non-K\"ahler almost Hermitian
structures, which are harmonic sections.

\noindent 1. (\cite{W1}) The standard nearly K\"ahler structure on
$S^6$.

\noindent 2. (\cite{W1}) The Calabi-Eckmann complex structure on
$S^{2p+1}\times S^{2q+1}$ with the product metric.

\noindent 3. (\cite{W1}) The Abbena-Thurston \cite{Ab,Tur} almost
K\"ahler structure on (the real Heisenberg group $\times
S^1$)/(discrete subgroup).

\noindent 4. The complex structure of the Iwasawa manifold - (the
complex Heisenberg group)/(discrete subgroup).

\smallskip

G. Bor, L. Hern\'andez-Lamoneda and M. Salvai \cite{BLS} have given
sufficient conditions for an almost Hermitian complex structure to
minimize the energy functional among sections of the twistor bundle.

\begin{tw}\rm{(\cite{BLS})}
Let $(N,h)$ be a compact Riemannian manifold and let $J$ be a compatible
almost complex structure on it. Suppose that

\noindent (1) $dim\,N=4$, the manifold $(N,h)$ is anti-self-dual and
the almost Hermitian structure $(h,J)$ is Hermitian, or almost
K\"ahler,

\noindent or

\noindent (2) $dim\,N\geq 6$, $(N,h)$ is conformally flat and the
almost Hermitian structure $(h,J)$ is of Gray-Hervella \cite{GH}
type $W_1\oplus W_4$.

Then the almost complex structure $J$ is an energy minimizer.
\end{tw}

\noindent {\bf Examples} of non-K\"ahler minimizers of the energy
functional.

\noindent 1. (\cite{BLS}) C. LeBrun \cite{LeBr} has constructed
anti-self-dual Hermitian structures on the blow-ups $(S^3\times
S^1)\sharp\, n\overline{{\mathbb C}{\mathbb P}^2}$ of the Hopf
surface $S^3\times S^1$. Blow-ups do not affect the first Betti
number, so any blow up of the Hopf surface has Betti number one, and
hence it does not admit a K\"ahler metric.

\noindent 2. I. Kim \cite{Kim} has shown the existence of
anti-self-dual strictly almost K\"ahler structures on ${\mathbb
C}{\mathbb P}^2\sharp\, n\overline{{\mathbb C}{\mathbb P}^2}$,
$n\geq 11$, $(S^2\times \Sigma)\sharp\, n\overline{{\mathbb
C}{\mathbb P}^2}$, genus $\Sigma\geq 2$, $(S^2\times T^2)\sharp\,
n\overline{{\mathbb C}{\mathbb P}^2}$, $n\geq 6$.

\noindent 3. (\cite{BLS}) The standard Hermitian structure on the
Hopf manifold $S^{2p+1}\times S^1$ is conformally flat and locally
conformally K\"ahler, and hence  of Grey-Hervella class $W_4$.

\smallskip

Let $(N,h,J)$ be an almost Hermitian manifold and
$\Omega(X,Y)=h(JX,Y)$ its  fundamental $2$-form. Then the
Euler-Lagrange equation $[J,\nabla^{\ast}\,\nabla J]=0$ is
equivalent to the identity
\begin{equation}\label{E-L-2}
(\nabla^{\ast}\nabla\,\Omega)(X,Y)=(\nabla^{\ast}\nabla\,\Omega)(JX,JY),\quad X,Y\in TN.
\end{equation}
Note that for the rough Laplacian $\nabla^{\ast}\nabla$ we have $\nabla^{\ast}\nabla\,\Omega=-Trace\,\nabla^2\Omega$.

\smallskip

The following simple observation is useful in many cases. Let
$\widehat\Omega$ be the section of $\Lambda^2TN$ corresponding to
the the $2$-form $\Omega$ under the isomorphism $\Lambda^2TN\cong
\Lambda^2T^{\ast}N$ determined by the metric $g$ on $\Lambda^2TN$
defined by means of the metric $h$ on $TN$ via (\ref{g^2}). Thus,
$h(\widehat\Omega,X\wedge Y)=\Omega(X,Y)$, and if
$E_1,...,E_m,JE_1,...,JE_m$ is an orthonormal frame of $TN$, then
$$
\widehat\Omega=2\sum_{k=1}^m E_k\wedge JE_k.
$$
Denote by ${\cal R}(\Omega)$ the $2$-form corresponding to
${\mathcal R}(\widehat\Omega)$; we have $${\mathcal
R}(\Omega)(X,Y)=h({\mathcal R}(\widehat\Omega),X\wedge Y).$$

\begin{lemma}\rm(\cite{DM02})\label{E-L-3}
A compatible almost complex structure $J$ on a Riemannian manifold $(N,h)$ is a harmonic section if and only if
\begin{equation}\label{Lap}
\Delta\Omega(X,Y)-\Delta\Omega(JX,JY)={\mathcal R}(\Omega)(X,Y)-{\mathcal R}(\Omega)(JX,JY),\quad X,Y\in TN,
\end{equation}
where $\Delta$ is the Laplace-de Rham operator of $(N,h)$.
\end{lemma}

\noindent {\bf Proof}. By the  Weitzenb\"ock formula
$$
\begin{array}{c}
\Delta\Omega(X,Y)-(\nabla^{\ast}\nabla\,\Omega)(X,Y)\\[6pt]
=Trace\{Z\to (R(Z,Y)\Omega)(Z,X)-(R(Z,X)\Omega)(Z,Y)\},
\end{array}
$$
$X,Y\in TN$ (see, for example, \cite{EL}). We have
$$
\begin{array}{c}
(R(Z,Y)\Omega)(Z,X)=-\Omega(R(Z,Y)Z,X)-\Omega(Z,R(Z,Y)X)\\[8pt]
=h(R(Z,Y)Z,JX)+h(R(Z,Y)X,JZ).
\end{array}
$$
Hence
$$
\begin{array}{c}
\Delta\Omega(X,Y)-(\nabla^{\ast}\nabla\,\Omega)(X,Y)\\[6pt]
=Ricci(Y,JX)-Ricci(X,JY)\\[6pt]
+Trace\,\{Z\to h(R(Z,Y)X,JZ-h(R(Z,X)Y,JZ)\}
\end{array}
$$
By the algebraic Bianchi identity
$$
\begin{array}{c}
h(R(Z,Y)X,JZ)-h(R(Z,X)Y,JZ)=h(R(X,Y)Z,JZ)\\[6pt]
=h({\mathcal R}(Z\wedge JZ),X\wedge Y).
\end{array}
$$
We have
$$
\begin{array}{c}
Trace\,\{Z\to h({\mathcal R}(Z\wedge JZ),X\wedge Y)\}=2\sum_{k=1}^mh({\mathcal R}(E_k\wedge JE_k),X\wedge Y)\\[6pt]
=h({\mathcal R}(\widehat\Omega),X\wedge Y)={\mathcal R}(\Omega)(X,Y).
\end{array}
$$
Thus
$$
\Delta\Omega(X,Y)-(\nabla^{\ast}\nabla\,\Omega)(X,Y)=Ricci(Y,JX)-Ricci(X,JY)+{\mathcal R}(\Omega)(X,Y),
$$
and the result follows from (\ref{E-L-2}).

\section{The Atiyah-Hitchin-Singer and Eells-Salamon almost complex structures as harmonic sections}

Lemma~\ref{E-L-3} has been used to prove the following statement.

\begin{tw} \rm(\cite{DM02})\label{har-sec}
Let $(M,g)$ be an oriented Riemannian $4$-manifold and let
$({\mathcal Z},h_t)$ be its negative twistor space. Then:

\noindent $(i)$ The Atiyah-Hitchin-Singer almost-complex structure
${\mathcal J}_1$ on $({\mathcal Z},h_t)$ is a harmonic section if
and only if $(M,g)$ is a self-dual manifold.

\noindent $(ii)$ The Eells-Salamon almost-complex structure ${\mathcal
J}_2$ on $({\mathcal Z},h_t)$ is a harmonic section if and only if
$(M,g)$ is a self-dual manifold with constant scalar curvature.
\end{tw}

In order to apply Lemma~\ref{E-L-3}, one needs to compute the
Laplacian of the fundamental $2$-form
$\Omega_{k,t}(A,B)=h_t({\mathcal J}_kA,B), k=1,2$, of the
almost-Hermitian structure $(h_t,{\mathcal J}_k)$ on ${\mathcal Z}$.
A computation, involving coordinate-free formulas for the
differential and codifferential of $\Omega_{k,t}$  (\cite{M}), gives
the following expression for the Laplacian of $\Omega_{k,t}$ in
terms of the base manifold $(M,g)$.

\begin{lemma}\rm(\cite{DM02}) \label{L-Omega} Let $V$ be a vertical vector of ${\cal Z}$ at a point $\sigma$ and
$X,Y\in T_{\pi(\sigma)}M$. Then
\begin{equation}\label{3.2}
\Delta\Omega_{k,t}(X^h,Y^h)_{\sigma}=g(\frac{4\sigma}{t}+2(-1)^k{\cal R}(\sigma),X\wedge Y)
+tg(R(X\wedge Y)\sigma,R(\sigma)\sigma)
\end{equation}
and
\begin{equation}\label{3.3}
\Delta\Omega_{k,t}(V,X^h)_{\sigma}=(-1)^{k+1}tg(\delta{\cal R}(X),V)
-tg((\nabla_{X}{\cal R})(\sigma),\sigma\times V).
\end{equation}
\end{lemma}
To compute the curvature terms ${\mathcal R_{\mathcal Z}}(\Omega_{k,t})$ in (\ref{Lap}) one can use the following
coordinate-free formula for the curvature of the twistor space.
\begin{prop}\rm (\cite{DM91})\label{sec}
Let ${\mathcal Z}$ be the negative twistor space of an oriented
Riemannian $4$-manifold $(M,g)$ with curvature tensor $R$. Let
$E,F\in T_{\sigma}{\mathcal Z}$ and $X=\pi_{\ast}E$,
$Y=\pi_{\ast}F$, $V={\cal V}E$, $W={\cal V}F$ where ${\cal V}$ means
"the vertical part". Then
$$
\begin{array}{c}
h_t(R_{{\mathcal Z}}(E,F)E,F)=g(R(X,Y)X,Y)\\[6pt]
-tg((\nabla_{X}{\cal R})(X\wedge Y),\sigma\times
W)+tg((\nabla_{Y}{\cal R})(X\wedge Y),\sigma\times V)\\[6pt]
-3tg({\cal R}(\sigma),X\wedge Y)g(\sigma\times V,W)\\[6pt]
-t^2g(R(\sigma\times V)X,R(\sigma\times
W)Y)+\displaystyle{\frac{t^2}{4}}||R(\sigma\times W)X+R(\sigma\times
V)Y||^2\\[6pt]
-\displaystyle{\frac{3t}{4}}||R(X,Y)\sigma||^2+t(||V||^2||W||^2-g(V,W)^2),
\end{array}
$$
where the norm of the vertical vectors is taken with respect to the metric $g$ on $\Lambda^2_{-}TM$.
\end{prop}
Using this formula, the well-known expression of the Levi-Civita
curvature tensor by means of sectional curvatures (cf. e.g.
\cite[\S\, 3.6, p. 93, formula (15)]{GKM}), and the differential
Bianchi identity one gets the following.

\begin{cotmb}\label{RZ}
Let $\sigma\in{\cal Z}$, $X,Y,Z,T\in T_{\pi(\sigma)}M$, and
$U,V,W,W'\in{\cal V}_{\sigma}$. Then
$$
\begin{array}{c}
h_t(R_{\cal
Z}(X^h,Y^h)Z^h,T^h)_{\sigma}=g(R(X,Y)Z,T)\\[6pt]
-\displaystyle\frac{3t}{12}[2g(R(X,Y)\sigma,R(Z,T)\sigma)
-g(R(X,T)\sigma,R(Y,Z)\sigma)\\[6pt]
+g(R(X,Z)\sigma,R(Y,T)\sigma)].
\end{array}
$$
$$
\begin{array}{c}
h_t(R_{\cal
Z}(X^h,Y^h)Z^h,U)_{\sigma}=-\displaystyle{\frac{t}{2}g(\nabla_{Z}{\cal
R}(X\wedge Y),\sigma\times U)}.\\[8pt]
h_t(R_{\cal
Z}(X^h,U)Y^h,V)_{\sigma}=\displaystyle{\frac{t^2}{4}g(R(\sigma\times
V)X,R(\sigma\times U)Y)}\\[6pt]
+\displaystyle{\frac{t}{2}g({\cal R}(\sigma),X\wedge
Y)g(\sigma\times V,U)}.\\[8pt]
h_t(R_{\cal Z}(X^h,Y^h)U,V)_{\sigma}
=\displaystyle{\frac{t^2}{4}}[g(R(\sigma\times V)X,R(\sigma\times
U)Y)\\[6pt]
\hspace{6cm}-g(R(\sigma\times U)X,R(\sigma\times V)Y)]\\[6pt]
+tg({\cal R}(\sigma),X\wedge Y)g(\sigma\times V,U)\\[8pt]
h_t(R_{\cal Z}(X^h,U)V,W)=0,\quad h_t(R_{\cal Z}(U,V)W,W')=g(U,W)g(V,W')-g(U,W')g(V,W).
\end{array}
$$
\end{cotmb}
This implies

\begin{lemma}\rm(\cite{DM02})\label{R-Omega} Let $V, W$ be vertical vectors of ${\cal Z}$ at a point $\sigma$ and
$X,Y\in T_{\pi(\sigma)}M$. Then
\begin{equation}\label{3.14}
\begin{array}{lll}
{\cal R}_t(\Omega_{k,t})(X^h,Y^h)_{\sigma}&=&2[1+(-1)^{k+1}]g({\cal R}(\sigma),X\wedge Y)-tg(R(X\wedge Y)\sigma,R(\sigma)\sigma)\\[6pt]
\vspace{0.15cm}
 & &-\displaystyle{\frac{t}{2}}Trace\{Z\to g(R(X\wedge Z)\sigma,R(Y\wedge K_{\sigma}Z)\sigma)\}\\[6pt]
\vspace{0.15cm}
 & &-\displaystyle{\frac{t}{2}}(-1)^kTrace\{{\cal V}_{\sigma}\ni\tau\to g(R(\tau)X,R(\sigma\times\tau)Y)\},                          \end{array}
\end{equation}

\noindent
 where the latter trace is taken with respect to the metric $g$ on ${\cal V}_{\sigma}$,

\begin{equation}\label{3.15}
{\cal R}_t(\Omega_{k,t})(V,X^h)_{\sigma}=tg((\nabla_{X}{\cal R})(\sigma),\sigma\times V)
\end{equation}
\vspace{0.15cm}
and
$$
\begin{array}{lll}
{\cal R}_t(\Omega_{k,t})(V,W)_{\sigma}&=&2[(-1)^{k+1}+tg({\cal R}(\sigma),\sigma)]g(V,\sigma\times W) \\[6pt]
& &+\displaystyle{\frac{t^2}{2}} Trace\{Z\to g(R(\sigma\times V)K_{\sigma}Z,R(\sigma\times W)Z)\}.
\end{array}
$$
\end{lemma}

\smallskip

\noindent {\bf Proof of Theorem~\ref{har-sec}}. According to  Lemmas~\ref{E-L-3},
\ref{L-Omega} and \ref{R-Omega}, the almost complex structure
${\mathcal J}_n$ is a harmonic section if and only if the following two conditions are satisfied:
\begin{equation}\label{3.16}
\begin{array}{c}
4g({\cal R}(\sigma),X\wedge Y-K_{\sigma}X\wedge K_{\sigma}Y)=\\[6pt]

tTrace\{Z\to g(R(X\wedge Z)\sigma,R(Y\wedge K_{\sigma}Z)\sigma)
 - g(R(K_{\sigma}X\wedge Z)\sigma,R(K_{\sigma}Y\wedge K_{\sigma}Z)\sigma)\} \\[8pt]

+t(-1)^kTrace\{{\cal V}_{\sigma}\ni\tau\to g(R(\tau)X,R(\sigma\times\tau)Y)
-g(R(\tau)K_{\sigma}X,R(\sigma\times\tau)K_{\sigma}Y)\}
\end{array}
\end{equation}
and
\begin{equation}\label{3.17}
g(\delta{\cal R}(K_{\sigma}X),\sigma\times V)=(-1)^kg(\delta{\cal R}(X),V)
\end{equation}
for every $\sigma\in{\cal Z}, V\in{\cal V}_{\sigma}$ and $X,Y\in T_{\pi(\sigma)}M$.

We shall show that condition (\ref{3.16}) is equivalent to $(M,g)$ being a self-dual manifold.
Note first that (\ref{3.16}) holds  for every $X,Y\in T_{\pi(\sigma)}M$ if and only if it holds for every
$X,Y\in T_{\pi(\sigma)}M$ with $||X||=||Y||=1$ and $X\perp Y, K_{\sigma}Y$. For every such $X,Y$ there is a
unique $\tau\in {\mathcal V}_{\sigma}$, $||\tau||=1$, such that $Y=K_{\tau}X$, namely
$\tau=X\wedge Y-K_{\sigma}X\wedge K_{\sigma}Y$; conversely, if $\tau\in {\mathcal V}_{\sigma}$, $||\tau||=1$
and $Y=K_{\tau}X$, then $X\perp Y, K_{\sigma}Y$ in view of (\ref{KK}). Thus, (\ref{3.16}) holds if and only if
it holds for every $X\in T_{\pi(\sigma)}M$ and $Y=K_{\tau}X$ with $||X||=1$, $\tau\in {\mathcal V}_{\sigma}$, $||\tau||=1$.
Given such $X$ and $\tau$, the vectors $E_1=X$, $E_2=K_{\sigma}X$, $E_3=K_{\tau}X$, $E_4=K_{\sigma\times\tau}X$
constitute an oriented orthonormal basis of $T_{\pi(\sigma}M$ such that $s_1^{-}=\sigma$, $s_2^{-}=\tau$,
$s_3^{-}=\sigma\times\tau$, where $s_1^{-},s_2^{-},s_3^{-}$ are defined by means of $\{E_1,...,E_4\}$ via (\ref{s-basis}).
Using the bases $\{E_1,...,E_4\}$ of $T_{\pi(\sigma)}M$ and $\tau, \sigma\times\tau$ of ${\mathcal V}_{\sigma}$ to compute the
traces in the right-hand side of (\ref{3.16}),  we see that identity (\ref{3.16}) is equivalent to
$$
\begin{array}{c}
4g({\mathcal R}(\sigma),\tau)=tg(R(\sigma)\sigma,R(\tau)\sigma)\\[6pt]
+t(-1)^{k}g(R(\tau)\sigma,R(\sigma\times\tau)(\sigma\times\tau))
-t(-1)^{k}g(R(\tau)(\sigma\times\tau),R(\sigma\times\tau)\sigma)
\end{array}
$$
for every  $\sigma,\tau\in{\mathcal  Z}$, $\pi(\sigma)=\pi(\tau)$, $\sigma\perp\tau$.
Using (\ref{eq 2.5}) we easily see also that the latter identity is equivalent to
\begin{equation}\label{3.18}
\begin{array}{c}
4g({\mathcal R}(\sigma),\tau)=tg({\mathcal R}(\sigma),\sigma\times\tau)g({\mathcal R}(\tau),\sigma\times\tau)
+tg({\mathcal R}(\sigma),\tau)g({\mathcal R}(\tau),\tau)\\[6pt]
+t(-1)^{k+1}g({\mathcal R}(\tau),\sigma\times\tau)g({\mathcal R}(\sigma\times\tau),\tau)
-t(-1)^{k+1}g({\mathcal R}(\tau),\sigma)g({\mathcal R}(\sigma\times\tau),\sigma\times\tau).
\end{array}
\end{equation}
 Writing this identity with $(\sigma,\tau)$
replaced by $(\tau,\sigma)$ and comparing the obtained identity  with (\ref{3.18}) we get
\begin{equation}\label{3.19}
g({\cal R}(\sigma),\tau)[g({\cal R}(\sigma),\sigma)-g({\cal R}(\tau),\tau)]=0.
\end{equation}
Replacing the pair $(\sigma,\tau)$ by $\displaystyle{\bigg(\frac{3\sigma+4\tau}{5},\frac{4\sigma-3\tau}{5}\bigg)}$ in
(\ref{3.19}) and using again this identity, we obtain
$$
[g({\cal R}(\sigma),\sigma)-g({\cal R}(\tau),\tau)]^2=4[g({\cal R}(\sigma),\tau)]^2,
$$
which, together with (\ref{3.19}), gives
$$
g({\cal R}(\sigma),\sigma)=g({\cal R}(\tau),\tau), \quad g({\cal R}(\sigma),\tau)=0.
$$
Thus
$$
g({\cal W}_{-}(\sigma),\sigma)=g({\cal W}_{-}(\tau),\tau)~ \mbox{ and } ~ g({\cal W}_{-}(\sigma),\tau)=0
$$
Since $Trace \, {\cal W}_{-}=0$, this implies ${\mathcal W}_{-}=0$.

    Conversely, if ${\cal W}_{-}=0$ we have ${\mathcal R}(\sigma)=\displaystyle{\frac{s}{6}}\sigma+{\mathcal B}(\sigma)$ where
${\mathcal B}(\sigma)\in\Lambda^2_{+}TM$, so it is obvious that identity (\ref{3.18}) is satisfied.

 To analyze condition (\ref{3.17}) we recall that $\delta{\mathcal R}=2\delta{\mathcal B}\>(=-dRicci)$
(cf. e.g. \cite{Besse}), so it follows from (\ref{dec}) that
$$
\delta{\cal R}(X)=-\frac{1}{3}{\rm grad}\>s\wedge X+2\delta{\cal W}(X),\>X\in TM.
$$
Suppose ${\cal W}_{-}=0$. Since $\delta{\cal W}_{+}(X)\in\Lambda^2_{+}TM$, we have
$$
g(\delta{\cal R}(X),V)=\frac{1}{3}g(X\wedge {\rm grad}\>s,V)
$$
for any $V\in\Lambda^2_{-}TM$. The latter formula and (\ref{eq 2.6}) imply that condition (\ref{3.17}) is equivalent
(for self-dual manifolds) to the identity
$$
g(V,X\wedge grad\>s)=(-1)^{k+1}g(V,X\wedge {\rm grad}\>s).
$$
Clearly, this identity is satisfied if $k=1$; for $k=2$ it holds if and only if the scalar curvature $s$ is constant.

\section{The second fundamental form of an almost Hermitian structure as a map into the twistor space}

Let $J$ be a compatible almost complex structure on a Riemannian manifold $(N,h)$.
Then we have a map
$J:(N,h)\to ({\mathscr T},h_t)$ between Riemannian manifolds. Let
$J^{\ast}T{\mathscr T}\to N$ be the pull-back of the bundle $T{\mathscr T}\to {\mathscr T}$
under the map $J:N\to {\mathscr T}$. We can consider
the differential $J_{\ast}:TN\to T{\mathscr T}$ as a section of the
bundle $Hom(TN,J^{\ast}T{\mathscr T})\to N$. Denote by $\widetilde D$
the connection on $J^{\ast}T{\mathscr T}$ induced by the Levi-Civita
connection $D$ on $T{\mathscr T}$. The Levi Civita connection $\nabla$
on $TN$ and the connection $\widetilde D$ on $J^{\ast}T{\mathscr T}$
induce a connection $\widetilde\nabla$ on the bundle
$Hom(TN,J^{\ast}T{\mathscr T})$. Recall that the second fundamental form
of the map $J$ is, by definition, $\widetilde\nabla J_{\ast}$.
The map $J: (N,h)\to ({\mathscr T},h_t)$ is harmonic if and only if
$$
Trace_{h}\widetilde\nabla J_{\ast}=0.
$$
Recall also that the map $J: (N,h)\to ({\mathscr T},h_t)$ is totally
geodesic exactly when $\widetilde\nabla J_{\ast}=0$.

\begin{prop}\rm(\cite{DM16, DHM15})\label{covder-dif}
For every  $X,Y\in T_pN$,
$$
\begin{array}{c}
\widetilde\nabla J_{\ast}(X,Y) =\displaystyle{\frac{1}{2}}{\mathcal
V}(\nabla^{2}_{XY}J + \nabla^{2}_{YX}J)\\[8pt]
-\displaystyle{\frac{2t}{n}}[(R((J\circ
\nabla_{X}J)^{\wedge})Y)^h_{J(p)}
+(R((J\circ\nabla_{Y}J)^{\wedge})X)^h_{J(p)}],
\end{array}
$$
where ${\mathcal V}$ means "the vertical component", $n=dim\,N$, and $\nabla^{2}_{XY}J=\nabla_X\nabla_Y J-\nabla_{\nabla_XY}J$ is
the second covariant derivative of $J$.
\end{prop}
The computation of the second fundamental form is based on several lemmas.

\smallskip

First, we note that identity (\ref{eq 2.5}) can be generalized as follows.

\begin{lemma}\label{R[a,b]} {\rm (\cite{D})} For every $a,b\in A(T_pN)$ and $X,Y\in T_pN$, we
have
\begin{equation}\label{Rw}
G(R(X,Y)a,b)=\frac{2}{n}h(R([a,b]^{\wedge})X,Y).
\end{equation}
\end{lemma}

{\bf Proof}. Let $E_1,...,E_n$ be an orthonormal basis of $T_pN$.
Then
$$
[a,b]=\frac{1}{2}\sum_{i,j=1}^n h([a,b]E_i,E_j)E_i\wedge E_j.
$$
Therefore
$$
\begin{array}{c}
h(R([a,b]^{\wedge})X,Y)=\displaystyle{\frac{1}{2}}\sum_{i,j=1}^n
h(R(X,Y)E_i,E_j)[h(abE_i,E_j)+h(aE_i,bE_j)]\\[8pt]
=\displaystyle{-\frac{1}{2}\sum_{i=1}^n
h(a(R(X,Y)E_i),bE_i)+\frac{1}{2}\sum_{k=1}^n
h(R(X,Y)aE_k,bE_k)}\\[8pt]
=\displaystyle{\frac{n}{2}G(R(X,Y)a,b)}.
\end{array}
$$

Lemma~\ref{R[a,b]}  implies

\begin{equation}\label{r-r}
h_t(R(X,Y)J,V)=\frac{2t}{n}h(R([J,V]^{\wedge})X,Y)=
\frac{4t}{n}h(R((J\circ V)^{\wedge})X,Y).
\end{equation}

\begin{lemma} \rm (\cite{D,DM91})\label{LC}
If $X,Y$ are vector
fields on $N$, and $V$ is a vertical vector field on ${\mathscr T}$, then
\begin{equation}\label{D-hh}
(D_{X^h}Y^h)_{I}=(\nabla_{X}Y)^h_{I}+\frac{1}{2}R_{p}(X\wedge Y)I
\end{equation}
\begin{equation}\label{D-vh}
(D_{V}X^h)_{I}={\cal H}(D_{X^h}V)_{I}=-\frac{2t}{n}(R_{p}((I\circ
V_I)^{\wedge})X)_I^h \ ,
\end{equation}
where $I\in {\mathscr T}$, $p=\pi(I)$, $n=dim\,N$, and ${\mathcal H}$ means "the
horizontal component".
\end{lemma}

{\bf Proof}. Identity (\ref{D-hh}) follows from the Koszul formula
for the Levi-Civita connection and the identity
$[X^h,Y^h]_I=[X,Y]^h_I+R(X,Y)I$.

Let $W$ be a vertical vector field on ${\mathscr T}$. Then
$$
h_t(D_{V}X^h,W)=-h_t(X^h,D_{V}W)=0 ,
$$
since the fibres are totally geodesic submanifolds, so $D_{V}W$ is a
vertical vector field. Therefore, $D_{V}X^h$ is a horizontal vector
field. Moreover, $[V,X^h]$ is a vertical vector field, hence
$D_{V}X^h={\cal H}D_{X^h}V$. Thus
$$
h_t(D_{V}X^h,Y^h)=h_t(D_{X^h}V,Y^h)=-h_t(V,D_{X^h}Y^h).
$$
Now (\ref{D-vh}) follows from (\ref{D-hh}) and (\ref{r-r}).

Any (local) section $a$
of the bundle $A(TN)$ determines a (local) vertical vector field
$\widetilde a$ on ${\mathscr T}$ defined by
$$
{\widetilde a}_I=\frac{1}{2}(a(p)+I\circ a(p)\circ I),\quad
p=\pi(I).
$$

The next lemma is "folklore".

\begin{lemma}\label{Xh-a til}
If $I\in{\mathscr T}$ and $X$ is a vector field on a neighbourhood of
the point $p=\pi(I)$, then
$$
[X^h,\widetilde a]_I=(\widetilde{\nabla_{X}a})_I.
$$
\end{lemma}

Let $I\in{\mathscr T}$ and let $U,V\in{\cal V}_I$. Take section $a$ and
$b$ of $A(TN)$ such that $a(p)=U$, $b(p)=V$ for $p=\pi(I)$. Let
$\widetilde a$ and $\widetilde b$ be the vertical vector fields
determined by the sections $a$ and $b$. Taking into account the fact
that the fibre of ${\mathscr T}$ through the point $I$ is a totally
geodesic submanifold, one easily gets by means of the Koszul formula that
\begin{equation}\label{a-b-til}
(D_{\widetilde a}\widetilde b)_I=\frac{1}{4}[UVI+IVU+I(UVI+IVU)I]=0.
\end{equation}

\begin{lemma}
\label{v-frame} For every $p\in N$, there exists an $h_t$-orthonormal
frame of vertical vector fields
$\{V_{\alpha}:~\alpha=1,...,m^2-m\}$, $m=\frac{1}{2}dim\,N$, in a
neighbourhood of the point $J(p)$ such that

\noindent $(1)$ $\quad (D_{V_{\alpha}}V_{\beta})_{J(p)}=0$,~~
$\alpha,\beta=1,...,m^2-m$.

\noindent $(2)$ $\quad$ If $X$ is a vector field near the point $p$,
then $[X^h,V_{\alpha}]_{J(p)}=0$.

\noindent $(3)$ $\quad$  $\nabla_{X_p}(V_{\alpha}\circ J)\perp {\cal
V}_{J(p)}.$

\end{lemma}

{\bf Proof}. Let $E_1,...,E_n$ be an orthonormal frame of $TN$ in a
neighbourhood  of $p$ such that $J(E_{2k-1})_p=(E_{2k})_p$,
$k=1,...,m$, and $\nabla E_l|_p=0$, $l=1,...,n$. Define sections
$S_{ij}, 1\leq i,j\leq n$, of $A(TN)$ by the formula
$$
S_{ij}E_l=\sqrt{\frac{n}{2}}(\delta_{il}E_j - \delta_{lj}E_i),\quad
l=1,...,n.
$$
Then $S_{ij}, i<j,$ form an orthonormal frame of $A(TN)$
with respect to the metric
$G(a,b)=\displaystyle{-\frac{1}{n}}Trace\,(a\circ b)\, ;a,b\in
A(TN)$. Set
$$
\begin{array}{c}
A_{r,s}=\frac{1}{\sqrt 2}(S_{2r-1,2s-1}-S_{2r,2s}),\quad
B_{r,s}=\frac{1}{\sqrt 2}(S_{2r-1,2s}+S_{2r,2s-1}),\\[6pt]
r=1,...,m-1,\>s=r+1,...,m.
\end{array}
$$
Then $\{(A_{r,s})_p,(B_{r,s})_p\}$ is a $G$-orthonormal basis of the
vertical space ${\cal V}_{J(p)}$. Note also that $\nabla
A_{r,s}|_p=\nabla B_{r,s}|_p=0$. Let $\widetilde A_{r,s}$ and
$\widetilde B_{r,s}$ be the vertical vector fields on ${\mathscr T}$
determined by the sections $A_{r,s}$  and $B_{r,s}$ of $A(TN)$.
These vector fields constitute a frame of the vertical bundle ${\cal
V}$ in a neighbourhood of the point $J(p)$.

Considering $\widetilde A_{r,s}\circ J$ as a section of $A(TN)$,  we have

$$
\begin{array}{c}
\nabla_{X_p}(\widetilde A_{r,s}\circ
J)=\frac{1}{2}\{(\nabla_{X_p}J)\circ
(A_{r,s})_p\circ J_p+J_p\circ (A_{r,s})\circ (\nabla_{X_p}J)\}\\[6pt]
=\frac{1}{2}\{-\nabla_{X_p}\circ J_p\circ (A_{r,s})_p+J_p\circ (A_{r,s})\circ (\nabla_{X_p}J)\}\\[6pt]
=\frac{1}{2}[(B_{r,s})_p,\nabla_{X_p}J].
\end{array}
$$
For every $I\in {\mathscr T}$, we have the orthogonal decomposition
\begin{equation}\label{ordec}
A(T_{\pi(I)}N)={\mathcal V}_I\oplus\{S\in A(T_{\pi(I)}N):~IS-SI=0\}.
\end{equation}
The endomorphisms $(B_{r,s})_p$ and $\nabla_{X_p}J$ of $T_pN$ belong
to ${\cal V}_{J(p)}$, so they anti-commute with $J(p)$, hence their commutator commutes with $J(p)$.
Therefore the commutator $[(B_{r,s})_p,\nabla_{X_p}J]$ is $G$-orthogonal to
the vertical space at $J$. Thus
$$
\nabla_{X_p}(\widetilde A_{r,s}\circ J)\perp{\cal V}_{J(p)},
$$
and, similarly,  $\nabla_{X_p}(\widetilde B_{r,s}\circ J)\perp{\cal
V}_{J(p)}$.

It is convenient to denote the elements of the frame $\{\widetilde
A_{r,s},\widetilde B_{r,s}\}$ by $\{\widetilde V_1,...,\widetilde
V_{m^2-m}\}$. In this way we have a frame of vertical vector fields
near the point $J(p)$ with  property $(3)$ of the lemma. Properties
$(1)$ and $(2)$ are also satisfied by this frame according to
(\ref{a-b-til}) and Lemma~\ref{Xh-a til}, respectively. In
particular,
$$
(\widetilde V_{\gamma})_{J(p)}(h_t(\widetilde V_{\alpha},\widetilde
V_{\beta}))=0,\quad \alpha,\beta,\gamma=1,...,m^2-m.
$$

Note also that, in view of (\ref{D-vh}),
$$
{\cal V}(D_{X^h}\widetilde V_{\alpha})_{J(p)}=[X^h,\widetilde
V_{\alpha}]_{J(p)}=0,
$$
hence
$$
X^h_{J(p)}(h_t(\widetilde V_{\alpha},\widetilde V_{\beta}))=0.
$$

 Now it is clear that the $h_t$-orthonormal frame
$\{V_1,...,V_{m^2-m}\}$ obtained from $\{\widetilde
V_1,...,\widetilde V_{m^2-m}\}$ by the Gram-Schmidt process has the
properties stated in the lemma.

\smallskip

\noindent {\bf Proof of Proposition~\ref{covder-dif}}.

Extend the tangent vectors $X$ and $Y$ to vector fields in a
neighbourhood of the point $p$. Let $V_1,...,V_{m^2-m}$ be an
$h_t$-orthonormal frame of vertical vector fields with properties
$(1)$ - $(3)$ stated in Lemma~\ref{v-frame}.

 We have
$$
J_{\ast}\circ Y=Y^h\circ J+\nabla_{Y}J=Y^h\circ
J+\sum_{\alpha=1}^{m^2-m}h_t(\nabla_{Y}J,V_\alpha\circ
J)(V_{\alpha}\circ J),
$$
hence
$$
\begin{array}{c}
{\widetilde D}_{X}(J_{\ast}\circ Y)=(D_{J_{\ast}X}Y^h)\circ J+
\sum_{\alpha=1}^{m^2-m}h_t(\nabla_{Y}J,V_{\alpha})(D_{J_{\ast}X}V_{\alpha})\circ J\\[8pt]
+t\sum_{\alpha=1}^{m^2-m}G(\nabla_X\nabla_{Y} J,V_{\alpha}\circ
J)(V_{\alpha}\circ J).
\end{array}
$$
This, in view of Lemma~\ref{LC}, implies
$$
\begin{array}{c}
\widetilde{D}_{X_p}(J_{\ast}\circ
Y)=(\nabla_{X}Y)^h_{J(p)}+\displaystyle{\frac{1}{2}}R(X\wedge
Y)J(p)-\displaystyle{\frac{2t}{n}}(R((J\circ\nabla_{X}J)^{\wedge})Y)^h_{J(p)}\\[8pt]
+t\sum_{\alpha=1}^{m^2-m}G(\nabla_{X_p}\nabla_{Y}J,V_{\alpha}\circ
J)_{p}V_{\alpha}(J(p))\\[8pt]
-\displaystyle{\frac{2t}{n}}(R((J\circ\nabla_{Y}J)^{\wedge})X)^h_{J(p)}\\[8pt]
=(\nabla_{X_p}Y)^h_{J(p)} +\displaystyle{\frac{1}{2}}{\cal
V}(\nabla_{X_p}\nabla_{Y}J +\nabla_{Y_p}\nabla_{X}J)
+\frac{1}{2}\nabla_{[X,Y]_p}J \\[8pt]
-\displaystyle{\frac{2t}{n}}[R((J\circ\nabla_{X}J)^{\wedge})Y)^h_{J(p)}
+(R((J\circ\nabla_{Y}J)^{\wedge})X)^h_{J(p)}].
\end{array}
$$
It follows that
$$
\begin{array}{c}
\widetilde\nabla J_{\ast}(X,Y)=\widetilde D_{X_p}(J_{\ast}\circ Y)-
(\nabla_{X}Y)^h_{\sigma}-\nabla_{\nabla_{X_p}Y} J\\[8pt]
=\displaystyle{\frac{1}{2}}{\cal
V}(\nabla_{X_p}\nabla_{Y}J-\nabla_{\nabla_{X_p}Y} J
+\nabla_{Y_p}\nabla_{X}J-\nabla_{\nabla_{Y_p}X} J)\\[8pt]
-\displaystyle{\frac{2t}{n}}[R((J\circ\nabla_{X}J)^{\wedge})Y)^h_{J(p)}
+(R((J\circ\nabla_{Y}J)^{\wedge})X)^h_{J(p)}].
\end{array}
$$

Proposition~\ref{covder-dif} implies immediately the following.
\begin{cotmb}
If $(N,h,J)$ is K\"ahler,  the map $J:(N,h)\to ({\mathscr T},h_t)$
is a totally geodesic isometric imbedding.
\end{cotmb}

\smallskip

\noindent {\bf Remark}. In view of the decomposition (\ref{ordec}),
the Euler-Lagrange equation $[J,\nabla^{\ast}\,\nabla J]=0$ is
equivalent to the condition that the vertical part of
$\nabla^{\ast}\nabla J=-Trace\,\nabla^2 J$ vanishes. Thus, by
Proposition~\ref{covder-dif}, $ J$ is a harmonic section if and only
if
$$
{\cal V}\,Trace\,\widetilde\nabla J_{\ast}=0.
$$

This fact, Proposition~\ref{covder-dif} and Theorem~\ref{har-sec}
imply
\begin{cotmb}\label{hs}

\noindent $(i)$ ${\cal V}\,Trace\,\widetilde\nabla {\cal
J}_{1\,\ast}=0$ if and only if $(M,g)$ is self-dual.

\noindent $(ii)$ ${\cal V}\,Trace\,\widetilde\nabla {\cal
J}_{2\,\ast}=0$ if and only if $(M,g)$ is self-dual and with
constant scalar curvature.
\end{cotmb}

\section{The Atiyah-Hitchin-Singer and Eells-Salamon almost complex structures as harmonic maps }

The main result in this section is the following.

\begin{tw}\label{H}
Each of the Atiyah-Hitchin-Singer  and
Eells-Salamon almost complex structures on the negative twistor space ${\cal Z}$ of an
oriented Riemannian  four-manifold $(M,g)$ determines a harmonic map
if and only if  $(M,g)$ is either self-dual and  Einstein, or
is locally the product of an open interval in ${\mathbb R}$ and a
$3$-dimensional Riemannian manifold of constant curvature.
\end{tw}

\medskip

\noindent{\bf Remarks}. 1. Every manifold that  is locally the product
of an open interval in ${\mathbb R}$ and a $3$-dimensional
Riemannian manifold of constant curvature $c$ is locally conformally
flat with constant scalar curvature $6c$. It is not Einstein unless
$c=0$, i.e. Ricci flat.

\noindent 2. According to Theorems 3 and 4, the conditions under which  ${\mathcal J_1}$ or ${\mathcal J_2}$  is  a harmonic section or a harmonic map do not depend on the parameter $t$ of the metric $h_t$. Taking certain special values of $t$, we can obtain a metric $h_t$ with nice properties (cf., for example, \cite{DM91, DGM, M}).

\medskip

The proof is based on several technical lemmas.

\smallskip

Note first that the almost complex structure ${\cal J}_k$, $k=1$ or
$2$, is a harmonic map if and only if ${\cal
V}\,Trace\,\widetilde\nabla {\cal J}_{k\,\ast}=0$ and ${\cal
H}\,Trace\,\widetilde\nabla {\cal J}_{k\,\ast}=0$. According to
Proposition~\ref{covder-dif},  ${\cal H}\,Trace\,\widetilde\nabla
{\cal J}_{k\,\ast}=0$, $k=1,2$, if and only if for every
$\sigma\in{\cal Z}$ and every $F\in T_{\sigma}{\cal Z}$
$$
Trace_{h_t}\,\{T_{\sigma}{\cal Z}\ni A\to h_t(R_{\cal Z}(({\cal
J}_{k}\circ D_{A}{\cal J}_k)^{\wedge})A),F)\}=0.
$$

Set for brevity
$$
Tr_k(F)=Trace_{h_t}\,\{T_{\sigma}{\cal Z}\ni A\to h_t(R_{\cal
Z}(({\cal J}_{k}\circ D_{A}{\cal J}_k)^{\wedge})A),F)\}.
$$

\smallskip

Let $\Omega_{k,t}(A,B)=h_t({\cal J}_kA,B)$ be the fundamental
$2$-form of the almost Hermitian manifold $({\cal Z},h_t,{\cal
J}_k)$, $k=1,2$. Then, for $A,B,C\in T_{\sigma}{\cal Z}$,
$$
h_t({\cal J}_{k}\circ D_{A}{\cal J}_k)^{\wedge},B\wedge
C)=-\frac{1}{2}h_t((D_{A}{\cal J}_k)(B),{\cal
J}_kC)=-\frac{1}{2}(D_{A}\Omega_{k,t})(B,{\cal J}_kC).
$$
\begin{lemma}\label{D-Omega}{\rm (\cite{M})}
Let $\sigma\in{\cal Z}$ and $X,Y\in T_{\pi(\sigma)}M$, $V\in{\cal
V}_{\sigma}$. Then
$$
(D_{X^h_{\sigma}}\Omega_{k,t})(Y^h_{\sigma},V)=\frac{t}{2}[(-1)^kg({\cal
R}(V),X\wedge Y)-g({\cal R}(\sigma\times V),X\wedge K_{\sigma}Y)],
$$
$$
(D_{V}\Omega_{k,t})(X^h_{\sigma},Y^h_{\sigma})=\frac{t}{2}g({\cal
R}(\sigma\times V),X\wedge K_{\sigma}Y +K_{\sigma}X\wedge
Y)+2g(V,X\wedge Y).
$$
Moreover, $(D_{A}\Omega_{k,t})(B,C)=0$ when $A,B,C$ are three
horizontal vectors at $\sigma$ or at least two of them are vertical.
\end{lemma}

\begin{cotmb}\label{JDJ}
Let $\sigma\in{\cal Z}$, $X\in T_{\pi(\sigma)}M$, $U\in {\cal
V}_{\sigma}$. If $E_1,...,E_4$ is an orthonormal basis of
$T_{\pi(\sigma)}M$ and $V_1,V_2$ is an $h_t$-orthonormal basis of
${\cal V}_{\sigma}$,
$$
\begin{array}{l}
({\cal J}_k\circ D_{X^h_{\sigma}}{\cal
J}_k)^{\wedge}=-\displaystyle{\frac{1}{2}}\sum_{i=1}^4\sum_{l=1}^2[g({\cal
R}(\sigma\times V_l),X\wedge E_i)\\[6pt]
\hspace{5.5cm} +(-1)^k g({\cal R}(V_l),X\wedge
K_{\sigma}E_i)](E_i^h)_{\sigma}\wedge V_l,
\end{array}
$$
$$
\begin{array}{l}
({\cal J}_k\circ D_{U}{\cal J}_k)^{\wedge}=\sum\limits_{1\leq i<j\leq
4}[\displaystyle{\frac{t}{2}}g({\cal R}(\sigma\times U),E_i\wedge
E_j-K_{\sigma}E_i\wedge K_{\sigma}E_j)\\[8pt]
\hspace{6cm}-2g(U,E_i\wedge K_{\sigma}E_j)](E_i^h)_{\sigma}\wedge
(E_j^h)_{\sigma}.
\end{array}
$$
\end{cotmb}

By Corollary~\ref{hs}, if the vertical part of
$Trace\,\widetilde\nabla {\cal J}_{k\,\ast}$ vanishes, then the
manifold $(M,g)$ is self dual. In the case when the base manifold is
self-dual, simple, but long computations involving
Corollary~\ref{JDJ}, the algebraic Bianchi identity, Corollary~\ref
{RZ},  and formula (\ref{eq 2.6}), give the next two lemmas, which
play an essential role in the proof of Theorem~\ref{H}.

\begin{lemma}\label{tr-ver}
Suppose that $(M,g)$ is self-dual. Then, if $U\in{\cal V}_{\sigma}$,
$$
Tr_k(U)=\frac{t}{4}g({\cal B}(U),{\cal B}(\sigma)), \> k=1,2.
$$
\end{lemma}

\begin{lemma}\label{tr-horr}
Suppose that $(M,g)$ is self-dual. Then, if $X\in T_{p}M$,
$p=\pi(\sigma)$,
$$
\begin{array}{c}
Tr_k(X^h_{\sigma})=[1+(-1)^k]\displaystyle{\frac{s(p)}{144}}X(s)
+\displaystyle{\frac{1}{12}}(\frac{ts(p)}{6}-2)X(s)\\[6pt]
+Trace_{h_t}\,\{{\cal V}_{\sigma}\ni V\to
[\displaystyle{\frac{t}{8}}g((\nabla_{X}{\cal B})(
V),{\cal B}( V))\\[6pt]
\hspace{7cm}+(-1)^{k+1}\displaystyle{\frac{ts(p)}{24}}g(\delta{\cal
B}(K_VX),V)]\}.
\end{array}
$$
\end{lemma}

\noindent {\bf Sketch of the proof of Theorem~\ref{H}}. Suppose that ${\cal J}_1$ or ${\cal J}_2$ is a harmonic
map. By Corollary~\ref{hs}, $(M,g)$ is self-dual or self-dual with
constant  scalar curvature. Moreover, $Tr_k(U)=0$ for every vertical
vector $U$,  and $Tr_k(X^h)=0$ for every horizontal vector $X^h$,
$k=1$ or $k=2$. Note that in both cases the first term in the
expression for $Tr_k(X^h)$  given in Lemma~\ref{tr-horr} vanishes.

Lemma~\ref{tr-ver} implies that  $||{\cal B}(\cdot)||^2=const$ on
every fibre ${\cal Z}_p$ of the twistor space. One can show that
this holds  if and only if, at every point $p\in M$, at least
three eigenvalues of the Ricci operator ${\rho}$ coincide. Then
the next step in the proof is to demonstrate that the condition
$Tr_k(X^h_{\sigma})=0$ for every $\sigma\in{\cal Z}$, $X\in
T_{\pi(\sigma)}M$, is equivalent to the pair of identities
\begin{equation}\label{Tr-B}
g(\delta{\cal B}(X),\sigma)=0
\end{equation}
and
\begin{equation}\label{psi}
\bigg(\frac{ts(p)}{144}+\frac{1}{6} \bigg)X(s)-\frac{t}{24}X(||\rho||^2 )=0.
\end{equation}
It is not hard to see that the identity (\ref{Tr-B}) is equivalent to
$$
Trace\,\{E\to g((\nabla_{E}\rho)(K_{\sigma}E),X)\}=0.
$$

Let $r(X,Y)$ be the Ricci tensor and set
$$
dr(X,Y,Z)=(\nabla_{Y}r)(Z,X)-(\nabla_{Z}r)(Y,X).
$$
Thus
$$
dr(X,Y,Z)=g((\nabla_{Y}\rho)(Z),X)-g((\nabla_{Z}\rho)(Y),X).
$$
Take an oriented orthonormal basis $(E_1,...,E_4)$ such that $E_2=K_{\sigma}E_1$ and
$E_4=-K_{\sigma}E_3$. Then
$$
dr(X,E_1,E_2)-dr(X,E_3,E_4)=\sum_{m=1}^4g((\nabla_{E_m}\rho)(K_{\sigma}E_m),X).
$$
Denote by $W_{-}$ the $4$-tensor corresponding to the operator
${\cal W}_{-}$,
$$
W_{-}(X,Y,Z,T)=g({\cal W}_{-}(X\wedge Y),Z\wedge T).
$$
By the differential Bianchi identity we have
\begin{equation}\label{dr}
dr(X,E_1,E_2)-dr(X,E_3,E_4)=-2[\delta W_{-}(X,E_1,E_2)-\delta
W_{-}(X,E_3,E_4)].
\end{equation}
Since $(M,g)$ is self-dual, we see from the latter identity that
identity (\ref{Tr-B}) is always satisfied. Identity (\ref{dr})
shows also that
\begin{equation}\label{dr=0}
dr(X,\sigma)=0,\quad \sigma\in{\cal Z},\quad X\in T_{\pi(\sigma)}M.
\end{equation}

Let $\lambda_1(p)\leq \lambda_2(p)\leq\lambda_3(p)\leq \lambda_4(p)$
be the eigenvalues of the symmetric operator $\rho_p:T_pM\to T_pM$
in ascending order. It is well-known that the functions
$\lambda_1,...,,\lambda_4$ are continuous (see, e.g. \cite[Chapter
Two, \S 5.7 ]{K} or \cite[Chapter I, \S 3]{R}). We have seen that,
at every point of $M$, at least three eigenvalues of the operator
$\rho$ coincide. The set $U$ of points at which exactly three
eigenvalues coincide is open by the continuity of
$\lambda_1,...,\lambda_4$. For every $p\in U$, denote the simple
eigenvalue of $\rho_p$ by $\lambda(p)$ and the triple eigenvalue by
$\mu(p)$, so the spectrum of $\rho$ is $(\lambda,\mu,\mu,\mu)$ with
$\lambda(p)\neq \mu(p)$ for every $p\in U$. As is well-known, the
implicit function theorem implies that the function $\lambda$ is
smooth. Then the function $\mu=\frac{1}{3}(s-\lambda)$ is also
smooth. It is also well-known that, in a neighbourhood of every
point $p$ of $U$, there is a (smooth) unit vector field $E$ which is
an eigenvector of $\rho$ corresponding to $\lambda$. (for a proof
see \cite[Chapter 9, Theorem 7]{L}).  Let $\alpha$ be the dual
$1$-form to $E$, $\alpha(X)=g(E,X)$. Then
$$
r(X,Y)=(\lambda-\mu)\alpha(X)\alpha(Y)-\mu g(X,Y)
$$
in a neighbourhood of $p$.  Using this representation of the Ricci
tensor, identity $\delta r=-\frac{1}{2}ds$, and (\ref{dr=0}) one can
prove that the  scalar curvature $s$ is locally constant on $U$.
Then identity (\ref{psi}) implies that $||\rho||^2$ is locally
constant. Thus, in a neighbouhood of every point $p\in U$, we have
$\lambda+3\mu=a$ and $\lambda^2+3\mu^2=b^2$, where $a$ and $b$ are
some constants. It follows that
$\mu=12^{-1}(3a\pm\sqrt{12b^2-3a^2})$. Note that $12b^2-3a^2\neq 0$,
since otherwise we would have $\mu=\frac{1}{4}a$, hence
$\lambda=a-3\mu=\frac{1}{4}a=\mu$,  a contradiction. Since $\mu$ is
continuous, we see that $\mu$ is constant, hence $\lambda$ is also
constant. Then one can show that the $1$-form $\alpha$ is parallel.
It follows that the restriction of the Ricci tensor to $U$ is
parallel.

In the interior of the closed set $M\setminus U$ the eigenvalues of
the Ricci tensor coincide, hence the metric $g$ is Einstein on this
open set. Therefore, the scalar curvature $s$ is locally constant on
$Int\,(M\setminus U)$ and the Ricci tensor is parallel on it. Thus,
the Ricci tensor is parallel on the open set $U\cup Int\,(M\setminus
U)=M\setminus bU$, where $bU$ stands for the boundary of $U$. Since
$M\setminus bU$ is dense in $M$, it follows that  the Ricci tensor
is parallel on $M$. This implies that the eigenvalues
$\lambda_1\leq...\leq \lambda_4$ of the Ricci tensor are constant.
Thus, either $M$ is Einstein, or exactly three of the eigenvalues
coincide. In the second case the simple eigenvalue $\lambda=0$ by
\cite[Lemma 1]{DGM}. Therefore $M$ is locally the product of an
interval in ${\mathbb R}$ and a $3$-dimensional manifold of constant
curvature.

Conversely, suppose that $(M,g)$ is self-dual and Einstein, or
locally  is the product of an interval and a manifold of constant
curvature. Then at least three of the eigenvalues of the Ricci
tensor coincide which, as we have noted, implies that $||{\cal
B}(\cdot)||^2=const$ on every  fibre of ${\cal Z}$. It follows that
$g({\cal B}(\sigma),{\cal B}({\tau}))=0$ for every
$\sigma,\tau\in{\cal Z}$ with $g(\sigma,\tau)=0$. Therefore,
$T_k(U)=0$ for every vertical vector $U$, $k=1,2$, by
Lemma~\ref{tr-ver}. Moreover, $T_k(X^h)=0$ by Lemma~\ref{tr-horr},
since the scalar curvature is constant and $\nabla{\cal B}=0$.

\section{Almost hermitian structures on $4$-manifolds that are harmonic maps}

Let $(M,g)$ be a Riemannian $4$-manifold and $J$  a compatible
almost complex structure on it. Henceforth in this section, we shall
consider $M$ with the orientation induced by $J$, and the {\it
positive} twistor space ${\mathcal Z}_{+}$ will be denoted by
${\mathcal Z}$.

\smallskip
Denote the Ricci tensor of $(M,g)$ by $\rho$ and let $\rho^{\ast}$
be the $\ast$-Ricci tensor of the almost Hermitian manifold
$(M,g,J)$. Recall that the latter is defined by
$$
\rho^{\ast}(X,Y)=trace\{Z\to R(JZ,X)JY\}.
$$

Denote by $N$ the Nijenhuis
tensor of $J$
$$N(Y,Z)=-[Y,Z]+[JY,JZ]-J[Y,JZ]-J[JY,Z].$$ It is
well-known (and easy to check) that
\begin{equation}\label{nJ}
2g((\nabla_XJ)(Y),Z)=d\Omega(X,Y,Z)-d\Omega(X,JY,JZ)+g(N(Y,Z),JX),
\end{equation}
for all $X,Y,Z\in TM$.

\subsection{The case of integrable $J$}

Suppose that the almost complex structure $J$ is integrable.
Denote by $B$ the vector field on $M$ dual to the Lee form $\theta=-\delta\Omega\circ J$ with respect to
the metric $g$. Then (\ref{nJ}) and the identity
$d\Omega=\theta\wedge\Omega$ imply the following well-known
formula
\begin{equation}\label{nablaJ}
2(\nabla_XJ)(Y)=g(JX,Y)B-g(B,Y)JX+g(X,Y)JB-g(JB,Y)X.
\end{equation}
It follows that, considering $J$ as a section of the vector bundle $\Lambda^2_{+}TM$,
$$
\nabla_X J=\displaystyle{\frac{1}{2}}(JX\wedge B+X\wedge JB).
$$
Using this formula and Proposition~\ref{covder-dif}, one can prove the following.
\begin{tw}\rm(\cite{DHM15})\label{harm-int}
Suppose that the almost complex structure $J$ is integrable. Then
the map $J:(M,g)\to ({\cal Z},h_t)$ is harmonic if and
only $d\theta$ is a $(1,1)$-form and $\rho(X,B)=\rho^{\ast}(X,B)$
for every $X\in TM$.
\end{tw}

\begin{cotmb}\label{Wood}
The map $J:(M,g)\to ({\cal Z},h_t)$ defined by an
integrable almost Hermitian structure $J$ on $(M,g)$ is a harmonic
section if and only if the $2$-form $d\theta$ is of type $(1,1)$.
\end{cotmb}

\smallskip

\noindent {\bf Remark}. The $2$-form $d\theta$ of a
Hermitian surface $(M,g,J)$ is of type $(1,1)$ if and only if the
$\star$-Ricci tensor $\rho^{\ast}$ is symmetric.

\medskip

The map $J:M\to {\cal Z}$ is an imbedding and one can ask when this
imbedding is minimal, i.e. when $J(M)$ is a minimal submanifold of
$({\mathcal Z},h_t)$. Note that $J$ is  minimal exactly when it is a
harmonic map from $M$ endowed with the metric $J^{\ast}h_t$ into
$(\mathcal Z,h_t)$.

If $\Pi$ is the second fundamental form
 of the submanifold  $J(M)$, then $\Pi(J_{\ast}X,J_{\ast}Y)$ is the normal
component of $\widetilde\nabla J_{\ast}(X,Y)$. In particular
$J(M)$ is a minimal submanifold if and only if the normal
component of $Trace\,\widetilde\nabla J_{\ast}$ vanishes.

\begin{tw} \rm(\cite{DHM15}) \label{min-int}
Suppose that the almost complex structure $J$ is integrable. Then
the map $J:M\to ({\cal Z},h_t)$ is a minimal isometric
imbedding if and only if  $d\theta$ is a $(1,1)$ form and
$\rho(X,B)=\rho^{\ast}(X,B)$ for every $X\perp \{B,JB\}$.
\end{tw}

\subsection{The case of symplectic $J$}

Suppose that  $(M,g,J)$ is almost K\"ahler (symplectic).

Denote by $\Lambda^2_{0}TM$ the subbundle of $\Lambda^2_{+}TM$
orthogonal to $J$ (thus $\Lambda^2_{0}T_pM={\cal V}_{J(p)}$).  Under this notation we have the following.

\begin{tw} \rm(\cite{DHM15})\label{harm-sympl}
Let $(M,g,J)$ be an almost K\"ahler $4$-manifold. Then the map
$J:(M,g)\to ({\cal Z},h_t)$ is harmonic if and only if the
$\ast$-Ricci tensor $\rho^{\ast}$ is symmetric and
$$
Trace\,\{\Lambda^2_{0}TM\ni\tau\to R(\tau)(N(\tau))\}=0.
$$
\end{tw}
The proof makes use of the Weitzenb\"ock formula.

\begin{tw}\rm(\cite{DHM15})\label{min-sympl}
Let $(M,g,J)$ be an almost K\"ahler  four-manifold. Then the map
$J:M\to ({\cal Z},h_t)$ is a minimal isometric imbedding, if and
only if the $\star$-Ricci tensor $\rho^{\ast}$ is symmetric, and for
every $p\in M$
$$
Trace\,\{\Lambda^2_{0}T_pM\ni\tau\to R_p(\tau)(N(\tau))\}\in{\cal
N}_p.
$$
\end{tw}

\subsection {Examples}\rm(\cite{DHM15})

\noindent {\bf Primary Kodaira surfaces}. Every primary Kodaira
surface $M$ can be obtained in the following way
\cite[p.787]{Kodaira}. Let $\varphi_k(z,w)$ be the affine
transformations of ${\mathbb C}^2$ given by
$$\varphi_k(z,w) = (z+a_k,w+\overline{a}_kz+b_k),$$ where $a_k$, $b _k$, $k=1,2,3,4$, are complex
numbers such that
$$
a_1=a_2=0, \quad Im(a_3{\overline a}_4) =m b_1\neq 0,\quad b_2\neq 0
$$
for some integer $m>0$. They generate a group $G$ of transformations
acting freely and properly discontinuously on ${\mathbb C}^2$, and
$M$ is the  quotient space ${\mathbb C}^2/G$.

It is well-known that $M$ can also be described as the quotient of
${\mathbb C}^2$ endowed with a group structure by a discrete
subgroup $\Gamma$. The multiplication on ${\mathbb C}^2$ is defined
by
$$
(a,b).(z,w)=(z+a,w+\overline{a}z+b),\quad (a,b), (z,w)\in  {\mathbb
C}^2,
$$
and $\Gamma$ is the subgroup generated by $(a_k,b_k)$, $k=1,...,4$ (see, for example, \cite{Borc}).

A frame of $\Gamma$-left-invariant vector fields on ${\mathbb
C}^2\cong {\mathbb R}^4$ is given by
$$
A_1=\frac{\partial}{\partial x}-x\frac{\partial}{\partial
u}+y\frac{\partial}{\partial v},\quad A_2=\frac{\partial}{\partial
y}-y\frac{\partial}{\partial u}-x\frac{\partial}{\partial v},\quad
A_3=\frac{\partial}{\partial u},\quad A_4=\frac{\partial}{\partial
v},
$$
where $x+iy=z$, $u+iv=w$. Let $g$ be the left-invariant Riemannian
metric on $M\cong {\mathbb C}^2/\Gamma$ obtained from the metric on
${\mathbb C}^2$ for which the frame $A_1,...,A_4$ is orthonormal.

It is a result by Hasegawa \cite{Has} that every complex structure
on $M$ is induced by a left-invariant complex structure on ${\mathbb
C}^2$. It is not hard to see (\cite{M86,D14}) that a left-invariant
almost complex structure $J$ on ${\mathbb C}^2$ compatible with the
metric $g$ is integrable if and only if it is given by
$$
JA_1=\varepsilon_1 A_2,\quad JA_3=\varepsilon_2 A_4,\quad \varepsilon_1,\varepsilon_2=\pm 1.
$$
It is easy to check that, by Theorem~\ref{harm-int}, the map $J:(M,g)\to ({\mathcal Z},h_t)$ is harmonic.

It is also easy to give an explicit description of the twistor space $({\mathcal Z},h_t)$ (\cite{D14}), since
$\Lambda^2_{+}M$ admits a global orthonormal frame defined by
$$
\begin{array}{c}
s_1=\varepsilon_1 A_1\wedge A_2+\varepsilon_2A_3\wedge A_4,\quad
s_2=A_1\wedge A_3+\varepsilon_1\varepsilon_2 A_4\wedge A_2,\\[6pt]
s_3=\varepsilon_2 A_1\wedge A_4+\varepsilon_1 A_2\wedge A_3.
\end{array}
$$
Then we have a natural diffeomorphism $F: {\cal Z}\cong M\times
S^2$ defined by  $\sum_{k=1}^3 x_k s_k(p))\to (p,x_1,x_2,x_3)$ under which
$J$ becomes the section $p\to (p,1,0,0)$. Denote the pushforward of the metric $h_t$ under $F$ again by
$h_t$. For $x=(x_1,x_2,x_3)\in S^2$, set
$$
u_1(x)=\varepsilon_1\varepsilon_2 (-x_3,0,x_1),\quad u_2(x)=\varepsilon_2(x_2,-x_1,0),
u_3(x)=0,\quad u_4(x)=\varepsilon_1 (0,x_3,-x_2).
$$
The differential $F_{\ast}$ sends the horizontal lifts
$A_i^h$ $i=1,...,4$, at a point $\sigma=\sum_{k=1}^3 x_k
s_k(p)\in{\cal Z}$ to the vectors $A_i+u_i$ of $TM\oplus TS^2$.
Then, if $X,Y\in T_pM$ and $P,Q\in T_xS^2$,
$$
\begin{array}{c}
h_t(X+P,Y+Q)=g(X,Y)\\[6pt]
+t<P-\sum_{i=1}^4 g(X,A_i)u_i(x),Q-\sum_{j=1}^4
g(Y,A_j)u_j(x)>
\end{array}
$$
where $<.,.>$ is the standard metric of ${\mathbb R}^3$.

Now suppose again that $J$ is a left-invariant almost complex
structure on ${\mathbb C}^2$ compatible with the metric $g$.
Then the almost Hermitian structure $(g,J)$ is almost K\"ahler
(symplectic) if and only if $J$ is given by (\cite{M86,D14})
$$
\begin{array}{c}
JA_1=-\varepsilon_1\sin\varphi A_3+\varepsilon_1\varepsilon_2\cos\varphi A_4,\quad
JA_2=-\cos\varphi A_3-\varepsilon_2\sin\varphi A_4, \\[6pt]
JA_3=\varepsilon_1\sin\varphi A_1+\cos\varphi A_2,\quad JA_4=-\varepsilon_1\varepsilon_2\cos\varphi
A_1+\varepsilon_2\sin\varphi A_2, \\[6pt]
\varepsilon_1,\varepsilon_2=\pm 1, \> \varphi\in [0,2\pi).
\end{array}
$$
Suppose that $J$ is determined by these identities and set
$$
\begin{array}{c}
E_1=A_1,\quad E_2=-\varepsilon_1\sin\varphi A_3+\varepsilon_1\varepsilon_2\cos\varphi A_4,\\[6pt]
E_3=\cos\varphi A_3+\varepsilon_2\sin\varphi A_4,\quad E_4=A_2.
\end{array}
$$
Then $E_1,...,E_4$ is an orthonormal frame of $TM$ for which
$JE_1=E_2$ and $JE_3=E_4$. Define an orthonormal frame
$s_l=s_l^{+}$, $l=1,2,3$, of $\Lambda^2_{+}TM$ by means of
$E_1,...,E_4$ via (\ref{s-basis}). Computing $\rho^{\ast}(E_i,E_j)$
one can see that the $\ast$-Ricci tensor is symmetric. Also,
computing the curvature and the Nijenhuis tensor, we have
$$
Trace\,\{\Lambda^2_{0}TM\ni\tau\to
R(\tau)(N(\tau))\}=R(s_2)(N(s_2))+R(s_3)(N(s_3))=0.
$$
Thus, by Theorem~\ref{min-sympl},  $J$ defines a harmonic map.

As in the preceding case, it is easy to find an explicit description
of the twistor space ${\cal Z}$ of $M$ and the metric $h_t$
(\cite{D14}). The frame $\}s_1,s_2,s_3\}$ gives rise to an obvious
diffeomeorphism $F: {\cal Z}\cong M\times S^2$ under which $J$ becomes the map $p\to (p,1,0,0)$.
The differential $F_{\ast}$ of this diffeomorphism sends the horizontal lifts
$E_i^h$, $i=1,...,4$, to $E_i+u_i$ where
$$
\begin{array}{c}
u_1(x)=(x_3\varepsilon_1\varepsilon_2\cos\varphi, x_3\varepsilon_2\sin\varphi,-x_1\varepsilon_1\varepsilon_2\cos\varphi-x_2\varepsilon_2\sin\varphi),\\[6pt]
u_2(x)=(x_2\varepsilon_1\varepsilon_2\cos\varphi,
-x_1\varepsilon_1\varepsilon_2\cos\varphi, 0) ,\quad
u_3(x)=(x_2\varepsilon_2\sin\varphi,-x_1\varepsilon_2\sin\varphi,0)\\[6pt]
u_4(x)=(-x_3\varepsilon_2\sin\varphi,
x_3\varepsilon_1\varepsilon_2\cos\varphi,x_1\varepsilon_2\sin\varphi-x_2\varepsilon_1\varepsilon_2\cos\varphi).
\end{array}
$$
for $x=(x_1,x_2,x_3)\in S^2$. Then, if $X,Y\in T_pM$ and $P,Q\in T_xS^2$,
\begin{equation}\label{ht-sympl}
\begin{array}{c}
h_t(X+P,Y+Q)=g(X,Y)\\[6pt]
+t<P-\sum_{i=1}^4 g(X,E_i)u_i(x),Q-\sum_{j=1}^4
g(Y,E_j)u_j(x)>.
\end{array}
\end{equation}

\subsection{Four-dimensional Lie groups} By a result of A. Fino \cite{F} for every
left-invariant almost K\"ahler structure $(g,J)$ with $J$-invariant
Ricci tensor on a $4$-dimensional Lie group $M$ there exists an
orthonormal frame of left-invariant vector fields $E_1,...,E_4$ such
that
$$
JE_1=E_2,\quad JE_3=E_4
$$
and
$$
\begin{array}{c}
\displaystyle{[E_1,E_2]=0,\quad [E_1,E_3]=sE_1+\frac{s^2}{t}E_2,\quad [E_1,E_4]=\frac{s^2-t^2}{2t}E_1-sE_2}\\[6pt]
\displaystyle{[E_2,E_3]=-tE_1-sE_2,\quad [E_2,E_4]=-sE_1-\frac{s^2-t^2}{2t}E_2,\quad [E_3,E_4]=-\frac{s^2+t^2}{t}E_3}
\end{array}
$$
where $s$ and $t\neq 0$ are real numbers. Using this table one can compute the $\ast$-Ricci and  Nijenhius tensors.
The computation shows that $J$ defines a harmonic map by virtue of Theorem~\ref{min-sympl}.

\subsection{ Inoue surfaces of type $S^0$} Let us recall the construction of these surfaces (\cite{Inoue}).
Take a matrix  $A\in SL(3,\mathbb{Z})$ with a
real eigenvalue $\alpha > 1$ and two complex eigenvalues $\beta$ and
$\overline{\beta}$, $\beta\neq\overline{\beta}$. Choose eigenvectors
$(a_1,a_2,a_3)\in\mathbb{R}^3$ and $(b_1,b_2,b_3)\in \mathbb{C}^3$
of $A$ corresponding to $\alpha$ and $\beta$, respectively. Then the
vectors $(a_1,a_2,a_3), (b_1,b_2,b_3),
(\overline{b_1},\overline{b_2},\overline{b_3})$  are
$\mathbb{C}$-linearly independent. Denote the upper-half plane in
$\mathbb{C}$ by ${\bf H}$ and let $\Gamma$ be the group of
holomorphic automorphisms of ${\bf H}\times {\mathbb C}$ generated
by
$$g_o:(w,z)\to (\alpha w,\beta z), \quad g_i:(w,z)\to (w+a_i,z+b_i), \>i=1,2,3 .$$
The group $\Gamma$ acts on ${\bf H}\times{\mathbb C}$ freely and
properly discontinuously.  Then $M=({\bf H}\times {\mathbb
C})/\Gamma$ is a complex surface known as the Inoue surface of type
$S^0$. It has been shown by F. Tricerri \cite{Tr} that every such a
surface admits a locally conformal K\"ahler metric $g$ (cf. also
\cite{DO}) obtained from the $\Gamma$-invariant Hermitian metric
$$
g=\frac{1}{v^2}(du\otimes du+dv\otimes dv)+v(dx\otimes dx+dy\otimes
dy),\quad u+iv\in{\bf H}, \quad x+iy\in {\mathbb C}.
$$
on ${\bf H}\times {\mathbb C}$.

By Corollary~\ref{Wood}, $J:(M,g)\to ({\cal Z},h_t)$ is a harmonic
section. It is also a minimal isometric imbedding by
Theorem~\ref{min-int}. However,  $J$ is not a harmonic map according
to Theorem~\ref{harm-int}.


\begin{thebibliography}{1000}

\bibitem{Ab} E.~Abbena, {\it An example of an almost K\"ahler manifold which is not K\"ahlerian},
Boll. Un. Mat. Ital. A(6) {\bf 3} (1984), 383-392.


\bibitem{AHS} M.~F.~Atiyah, N.~J.~Hitchin, I.~M.~Singer, {\it Self-duality in four-dimensional Riemannian geometry},
Proc. Roy. Soc. London, Ser.A 362 (1978), 425-461.

\bibitem{Besse} A.~Besse , {\it Einstein manifolds}, Classics in Mathematics, Springer-Verlag,  2008.

\bibitem{BLS} G.~Bor, L.~Hern\'andez-Lamoneda, M.~Salvai, {\it
Orthogonal almost-complex structures of minimal energy}, Geom.
Dedicata {\bf 127} (2007), 75-85.

\bibitem{Borc} C.~Borcea, {\it Moduli for Kodaira surfaces}, Composition Math. {\bf 52} (1984),
373-380.

\bibitem{CG} E.~Calabi, H.~Gluck, {\it What are the best almost-complex structures on the 6-sphere?},
Proc. Sym. Pure Math. {\bf 54} (1993), part 2, 99-106.


\bibitem{D} J.~Davidov, {\it Einstein condition and twistor spaces of compatible partially
complex structures}, Diff. Geom. and its Appl. {\bf 22} (2005),
159-179

\bibitem{DM91} J.~Davidov, O.~Mushkarov, {\it On the Riemannian curvature of a twistor space},  Acta Math. Hungarica
{\bf 58} (1991), 319-332.

\bibitem{DM02} J.~Davidov, O.~Mushkarov, {\it Harmonic almost-complex structures on twistor spaces}, Israel J. Math.
{\bf 131} (2002), 319-332.

\bibitem{DHM15} J.~Davidov, A.~Ul~Haq, O.~Mushkarov, {\it Almost complex structures that are harmonic maps},
arXiv:1504.01610v2 [math.DG] 19 Aug 2015.

\bibitem{DM16} J.~Davidov, O.~Mushkarov, {\it Harmonicity of the Atiyah-Hitchin-Singer  and Eells-Salamon
almost complex structures}, in preparation.

\bibitem{DGM} J.~Davidov, G.~Grantcharov, O.~Mushkarov, {\it Twistorial
examples of $\ast$-Einstein manifolds}, Ann. Glob. Anal. Geom. {\bf
20} (2001), 103-115.


\bibitem{D14} J.~Davidov, {\it Twistorial construction of minimal hypersurfaces}, Inter. J. Geom. Methods in Modern Physics
{\bf 11}, No 6 (2014), 1459964.

\bibitem{DO} S.~Dragomir, L.~Ornea, {\it Locally conformal K\"ahler
geometry}, Progress in Math., v. 155, Birkh\"auser,
Boston-Basel-Berlin, 1998.


\bibitem{EL} J.~Eells, L.~Lemaire, {\it Selected topics in harmonic maps}, Cbms Regional Conference
Series in Mathematics, vol. {\bf 50}, AMS, Providernce, Rhode Island, 1983.

\bibitem{ES} J.~Eells, S.~Salamon, {\it Twistorial constructions of harmonic maps of surfaces into four-manifolds},
Ann. Scuola Norm. Sup. Pisa, ser.IV, 12 (1985), 589-640.

\bibitem{F} A.~Fino, {\it Almost K\"ahler $4$-dimensional Lie groups with $J$-invariant Ricci tensor},
Diff. Geom. Appl. {\bf 23} (2005), 26-37.


\bibitem{G} A.~Gray, {\it Minimal varieties and almost Hermitian
manifolds}, Michigan Math. J. {\bf 12} (1965), 273-287.

\bibitem{GH} A.~Gray, L.~M.~Hervella, {\it The sixteen classes of almost
Hermitian manifolds and their linear invariants}, Ann. Mat. Pure
Appl. {\bf 123} (1980), 288--294.

\bibitem{GKM} D.~Gromol, W.~Klingenberg, W.~Meyer, {\it Riemannsche
Geometrie in Grossen}, Lecture Notes in Mathematics vol. {\bf 55},
Springer-Verlag, 1968.


\bibitem{Has} K.~Hasegawa, {\it Complex and K\"ahler structures on compact solvmanifolds}, J.
Symplectic Geom. {\bf 3} (2005), 749-767.

\bibitem{Inoue} M.~Inoue, {\it On surfaces of class $VII_0$}, Invent. Math. {\bf 24} (1974), 269-310.


\bibitem{K} T.~Kato, {\it Perturbation theory for linear operators},
Springer-Verlag, Berlin-Heidelberg-New York, 1980.

\bibitem{Kim} I.~Kim, {\it Almost K\"ahler anti-self-dual metrics},
Ph.D. thesis, Stony Brook University, May 2014, available at
www.math.stonybrook.edu/alumni/2014-Inyong-Kim.pdf; see also
arXiv:1511.07656v1 [math.DG] 24 Nov 2015.

\bibitem{Kodaira} K.~Kodaira, {\it On the structure of compact complex analytic surfaces I}. Amer. J. Math. {\bf 86} (1964),
751-798.


\bibitem{L} P.~Lax, {\it Linear algebra and its applications},
John \& Sons, Inc., Hoboken, New Jersey, 2007.

\bibitem{LeBr} C.~LeBrun, {\it Anti-self-dual hermitian metrics on
blow-up Hopf surfaces}, Math. Ann. {\bf 289} (1991), 383-392.

\bibitem{M86} O.~Mu\v skarov, {\it Two remarks on Thurston's example}, In: Complex analysis and applications '85 (Varna, 1985),
 Publ. House Bulgar. Acad. Sci., Sofia, 1986, pp. 461-468.


\bibitem{M} O.~Mu\v skarov, {\it  Structures presque hermitienes sur espaces
twistoriels et leur types}, C. R. Acad. Sci. Paris S\'er.I Math.
{\bf 305} (1987), 307-309.


\bibitem{R} F.~Rellich, {\it Perturbation theory of eigenvalue
problems}, Notes on mathematics and its applications, Gordon and
Breach science publishers, New York-London-Paris, 1969.

\bibitem{ST} I.~M.~Singer, J.~A.~Thorpe, {\it The curvature of
$4$-dimensional Einstein spaces}, in papers in Honor of K. Kodaira,
Princeton University Press (Princeton), 1969, pp. 355-365.

\bibitem{Tr} F.~Tricerri, {\it Some example of locally conformal
K\"ahler manifolds}, Rend. Sem. Mat. Univ. Torino {\bf 40} (19982),
81-92.

\bibitem{Tur} W.~P.~Thurston, {\it Some simple examples of symplectic
manifolds}, Proc. Amer. Math. Soc. {\bf 55} (1976), 467-468.

\bibitem{V} J.~Vilms, {\it Totally geodesic maps}, J. Diff. Geom. 4 (1970), 73-79.

\bibitem{W1}C.~M.~Wood, {\it Instability of the nearly-K\"ahler six-sphere}, J. reine angew. Math.
{\bf 439} (1993), 205-212.


\bibitem{W2} C.~M.~Wood, {\it Harmonic almost-complex structures}, Compositio Mathematica {\bf 99}
(1995), 183-212.



\end{thebibliography}
\end{document}